\documentclass[12pt]{article}
\usepackage{theorem,amsfonts,amssymb}

\newtheorem{theorem}{Theorem}
\newtheorem{proposition}{Proposition}
\newtheorem{lemma}{Lemma}
\newtheorem{corollary}{Corollary}
\theorembodyfont{\upshape}
\newtheorem{definition}{Definition}
\newtheorem{remark}{Remark}
\newtheorem{proof}{Proof}

\newtheorem{acknowledgement}{Acknowledgement}

\newcommand{\bt}{\begin{theorem}}
\newcommand{\et}{\end{theorem}}
\newcommand{\bl}{\begin{lemma}}
\newcommand{\el}{\end{lemma}}
\newcommand{\bp}{\begin{proposition}}
\newcommand{\ep}{\end{proposition}}
\newcommand{\bd}{\begin{definition}}
\newcommand{\ed}{\end{definition}}
\newcommand{\br}{\begin{remark}}
\newcommand{\er}{\end{remark}}
\newcommand{\bc}{\begin{corollary}}
\newcommand{\ec}{\end{corollary}}
\newcommand{\bo}{\begin{proof}}
\newcommand{\eo}{\end{proof}}
\newcommand{\be}{\begin{enumerate}}
\newcommand{\ee}{\end{enumerate}}

\newcommand{\Z}{{\mathbb Z}}

\newcommand{\mH}{{\mathbb H}}
\newcommand{\mP}{{\mathbb P}}
\newcommand{\mF}{{\mathbb F}}
\newcommand{\mE}{{\mathbb E}}
\newcommand{\mG}{{\mathbb G}}

\newcommand{\N}{{\mathbb N}}
\newcommand{\mL}{{\mathbb L}}
\newcommand{\R}{{\mathbb R}}
\newcommand{\C}{{\mathbb C}}
\newcommand{\mT}{{\mathbb T}}
\newcommand{\cA}{{\cal A}}
\newcommand{\cB}{{\cal B}}
\newcommand{\cH}{{\cal H}}
\newcommand{\cL}{{\cal L}}

\title{Recurrence and ergodicity of random walks on linear groups and on 
homogeneous spaces}
\author{Y. Guivarc'h and C. R. E. Raja}
\date{ }

\begin{document}
\maketitle

\let\epsi=\epsilon
\let\vepsi=\varepsilon
\let\lam=\lambda
\let\Lam=\Lambda
\let\ap=\alpha
\let\vp=\varphi
\let\ra=\rightarrow
\let\Ra=\Rightarrow
\let\da=\downarrow
\let\Llra=\Longleftrightarrow
\let\Lla=\Longleftarrow
\let\lra=\longrightarrow
\let\Lra=\Longrightarrow
\let\ba=\beta
\let\ga=\gamma
\let\Ga=\Gamma
\let\Da=\Delta
\let\un=\upsilon
\let\ol =\overline
\let\oa =\omega
\let\Oa=\Omega
\let\ct =\cdot
\let\bs =\backslash

\begin{abstract}
We discuss recurrence and ergodicity properties of random walks and
associated skew products for large classes of locally compact groups
and homogeneous spaces. In particular we show that a closed subgroup
of a product of finitely many linear groups over local fields supports
a recurrent random walk
if and only if it has at most quadratic growth.  We give also a
detailed analysis of ergodicity properties for special classes of
random walks on homogeneous spaces.  The structure of closed
subgroups of linear groups over local fields and the properties of
group actions with respect to stationary measures play an important role
in the proofs.
\end{abstract}



\begin{section}{Introduction}

Let $G$ be a locally compact second countable group and $\mu$ be a
Borel probability measure on $G$.  We denote by $\mP$ the product
measure $\mP = \mu ^{\N }$ on $\Omega = G^{\N}$.  Let $Y_i(\oa )$ be
the $i$-th coordinate of $\oa \in \Omega$ in $G$ ($i\in \N$). Then the left
random walk on $G$ of law $\mu$, starting from $g \in G$ is the
sequence of $G$-valued random variables $X_n^g$ defined by $X_n^g
(\oa )= Y_n (\oa )\cdots Y_1(\oa ) g$, $X_0^g(\oa )=g$.  Given a
$G$-space $E$ and $x \in E$, we write $X_n(\oa ,x) =X_n^e(\oa )x$.
The sequence $X_n(\oa ,x)$ is called the random walk of law $\mu$ on
$E$, starting from $x$ and its properties will play an essential
role in the study of $X_n^g(\oa )$ if $E$ is chosen as a homogeneous
space of $G$.  We say that $X_n^g$ is {\it
recurrent} if for every neighborhood $W$ of the identity $e$ of $G$
$\mP$-a.e $X_n(\oa )= X_n^e (\oa )\in W$ infinitely often.
 If $X_n^g$ is not recurrent, then $\mP$-a.e $X_n^g(\oa )$ escapes to
infinity and the random walk $X_n^g$ is said to be transient.  We
can also define the right random walk $^gX_n=gY_1 \cdots Y_n$.  Also
recurrence of $^gX_n$ is defined in a similar fashion.  From the
above we see that recurrence of $^gX_n$ is equivalent to recurrence
of $X_n^g$. Hence, in this case we say that $\mu$ is recurrent.  We
denote by $\mu ^k$ the $k$-th convolution power of $\mu$.  Then
recurrence of $\mu$ is equivalent to the condition that $\sum _0
^\infty \mu ^k(U)=\infty$ for some neighborhood $U$ of $e$ in $G$.  We
denote by $G_\mu$ the smallest closed subgroup of $G$ containing the
support of $\mu$ and we say that $\mu$ is adapted if $G_\mu =G$.
Then the group $G$ is said to be {\it recurrent} if there exists an
adapted probability $\mu$ on $G$ such that $\mu$ is recurrent.

A stronger notion of recurrence where Haar measure enters explicitly
in the definition is H-recurrence.  We will say that $G$ is {\it
H-recurrent} if there exists a probability measure $\mu $ on $G$
such that for every Borel set $B$ with positive Haar measure
$X_n(\oa ) \in B$, $\mP$-a.e infinitely often.

It is a classical result that ${\mathbb R}^p \times {\mathbb Z}^q$
is recurrent if and only if $p+q\leq 2$ (cf. \cite{Re}).  Also a
countable abelian group is recurrent if and only if it has rank at
most two (see \cite{Du}).  These groups are also $H$-recurrent.  We
denote by $\lam _G$ a left Haar measure on $G$.  We recall that $G$
is said to have polynomial growth of degree at most $d\in \N$ if for
each compact neighborhood $W$ of $e\in G$ there exists a constant
$c$ such that $\lam _G(W^n) \leq cn^d$ for every $n \in \N$; when
$d=2$ we say that $G$ has at most quadratic growth.  We observe that
if $G$ has polynomial growth, then $G$ as well as its closed
subgroups is unimodular.

We give now some reference from previous works.  The idea of
relating growth and recurrence appeared in \cite{Ke} for the case of
countable groups: non-exponential growth of finitely generated recurrent 
groups was conjectured.   The
case of locally compact groups and especially real Lie groups, was
considered in \cite{GKR} where it was shown that recurrent groups
are amenable and unimodular. Probabilistic ideas and calculations
were developed in order to show transience or recurrence in various
situations, particularly for connected Lie groups and their
countable subgroups. For example, it was shown there that the group
$\mG _2$ of euclidean motion of $\R ^2$ is recurrent while its
universal cover as well as the $3$-dimensional Heisenberg groups
(continuous or discrete) is not. Also the affine group of the line
was shown to be transient. In \cite{GK} the following "quadratic
growth conjecture" was stated: $G$ is recurrent if and only if it
has polynomial growth of degree at most two. This conjecture has
been settled for various classes of groups: compactly generated
nilpotent groups (cf. \cite{GK}), connected groups (cf. \cite{Ba}),
finitely generated groups (cf. \cite{VSC}), quasi-transitive groups
of automorphisms of graphs (cf. \cite{Wo}), p-adic Lie groups and
totally disconnected groups of polynomial growth (cf. \cite{Ra2}).
Powerful analytical techniques from real analysis were developed in
\cite{VSC} which allows to obtain precise asymptotics for $\mu ^n$
if $\mu$ has nice density and $G$ is compactly generated unimodular.
In particular, a direct relation between polynomial growth and
recurrence of such a $\mu$ was obtained.  See also \cite{Fur},
\cite{Gu2} and \cite{Wo} for surveys on random walks and recurrence
properties.

Here one of our main results proves the validity of the conjecture
for the class of closed subgroups of products of full linear groups
over a finite number of local fields.  We show also the validity of
the above conjecture in the following situations:
\be

\item[(a)] $G$ is a real Lie group;

\item[(b)] $G$ is H-recurrent.

\ee

In fact, we show that $G$ is compactly generated and has quadratic
growth is equivalent to the following:  $G$ has a compact normal
subgroup $K$ and a finite index subgroup $G_1$ such that $K\subset
G_1$ and $G_1/K$ is isomorphic to a closed subgroup of the motion
group of the plane.  In the situations considered here, if $G$ is
compactly generated and totally disconnected we get that up to
finite index and to a compact normal subgroup, $G$ is isomorphic to
a subgroup of $\Z ^2$. Concerning the above conjecture, our main
result is the following.

\bt \label{mt} Assume $G$ is a closed group of a product of finitely
many linear groups over local fields: $G\subset \prod _{i\in I}
GL(d_i, {\mathbb F}_i)$.  If $G$ is recurrent, then $G$ has at most
quadratic growth and in addition if $G$ is compactly generated, then
$G$ contains a compact normal subgroup $K$ and a finite index
subgroup $H$ such that $K\subset H$ and $H/K$ is isomorphic to a
closed subgroup of the group $\mG _2$ of euclidean motions of the
plane. \et

To our knowledge this result is new even in case of one field
$\mathbb F$, in particular for the real field. For the proof we use
in particular the following previously-known results. The analysis
developed in \cite {Ba} and \cite{GKR} for recurrence in locally
compact groups and real Lie groups.  The structure of finitely
generated recurrent groups given in \cite{VSC} which is based on
\cite{Gr}.  The existence of invariant Radon measures for a class of
random walks in homogeneous spaces as follows from Theorem 5.1 of
\cite{L}; this result leads to the unimodularity of closed subgroups
of recurrent groups. The structural results for automorphisms of
totally disconnected groups in \cite{BW} and \cite{JR} allows us to
use efficiently the facts above in case of linear groups over
non-archimedean fields.

In this paper we give also a detailed analysis of recurrence and
ergodicity properties for some special classes of homogeneous spaces
and random walks which occur naturally in other studies.  Let $E=
G/H$ be such a $G$-space and $\mu \in M^1(G)$.  As shown
below and in contrast to the group case, recurrence properties are
valid for $X_n(\oa )x$ in a much wider setting not related to
polynomial growth of $G$ and the asymptotic properties of $X_n(\oa )x$
depend strongly on $x$,  in general.  Then it is convenient to discuss 
recurrence properties in terms of a fixed Radon measure 
$\lam$ such that the class of
$\lam$ is $\mu$-invariant. For a positive Radon measure $\eta$ on
$E$, we define the convolution of $\mu$ and $\eta$ by $\mu
*\eta = \int g\eta d\mu (g)$ where $g\eta$ is the push-forward of
$\eta$ by $g\in G$.  Let $\lam$ be a fixed $\mu$-invariant
measure on $E$ $(\mu*\lam = \lam$), $\hat \Oa = G^\Z = G^{-\N } \times
G \times \Oa$ be the product
space, $\theta$ be the shift on $\hat \Oa $and $\oa \in \Oa$ the projection
of $\hat \oa \in \hat \Oa$ on $\Oa $.  We consider the skew product
$(\Oa \times E, \tilde \theta , \mP \otimes \lam)$ where
$\tilde \theta $ is defined by $\tilde \theta (\oa , x) = (\theta
\oa , Y_1 (\oa ) x)$ for $(\oa , x)\in \Oa \times E$.  We consider also the map
$\hat \theta $ of $\hat \Oa \times E$ into itself defined by
$\hat \theta (\hat \oa ,x) = (\theta \hat \oa , Y_1(\oa )x)$ for
$(\hat \oa ,x) \in \hat \Oa \times E$.  The second coordinate of
$\hat \theta ^n (\hat \oa , x)$ for $n \in
\Z$ defines an extension of $X_n(\oa ,x)$ to negative time, i.e a
bilateral random walk on $E$.  We observe that $\tilde \lam =
\mP \otimes \lam$ is $\tilde \theta $-invariant and we will show
that there exists a unique $\hat \theta $-invariant measure
$\hat \lam$ on $\hat \Oa \times E$ which has projection $\tilde \lam$
on $\Oa \times E$.  The system $(\hat \Oa \times E, \hat \theta ,
\hat \lam)$ can be considered as the natural extension of the system
$(\Oa \times E, \tilde \theta , \tilde \lam)$.  In section 2.3 we study
recurrence and ergodicity properties of such systems from a general
point of view. As a result of this discussion in Section 5 we prove
the following

\bt\label{emt}
In the situation of examples 2, 4, 5 of Section 5, the skew product $(\hat
\Oa \times E, \hat \theta , \hat \lam )$ is ergodic with respect to
the infinite $\hat \theta$-invariant measure $\hat \lam$.
\et

This result provides large classes of invertible transformations
with stochastic properties which are ergodic with respect to natural
infinite invariant measures.  Examples 4 and 5 can be considered as fibered 
dynamical systems with fiber $\R$ or $\Z$, which have properties similar to those 
of the skew products considered in \cite{Gu3}.

As an illustration of recurrence and ergodicity properties on homogeneous
spaces we show the
singularity of stationary measure for the action of $SL(2, \Z)$ and its
cofinite subgroups on the projective line.

\begin{acknowledgement}
The authors wish to thank School of Mathematics, TIFR, Mumbai India
for its hospitality during the special trimester Lie Group-Ergodic
theory and Probability measures.  The authors would also like to
thank B. Bekka, J.P. Conze and G. Willis for useful remarks.  Also, the 
first named author would like to thank ESI (Vienna) for hospitality during the 
special semester on Random walks and amenability (2007).  
\end{acknowledgement}

\end{section}

\begin{section}{Recurrence of random walks on groups and $G$-spaces}

\begin{subsection}{Some basic facts}

A left-random walk $X_n^g$ on $G$ is said to be transient (resp.
recurrent) if for any compact neighborhood $W$ of $e \in G$ we have
$\mP$-a.e, $X_n ^e (\oa )\not \in W$ for all large $n$ (resp.
$X_n^e(\oa ) \in W$, infinitely often): see \cite{Fur}, \cite{Gu2}
\cite{Re} and \cite{Wo} for further detailed information. A random
walk is either transient or recurrent.  Transience of $X_n^g$ is
equivalent to the fact that $\sum _{0}^\infty \mu ^k$ is finite on
compact sets.

Given a Borel probability measure $\mu$ on $G$, a Markov operator
$P=P_\mu$ on $G$ is defined by $P_\mu \psi (x) = \int \psi (gx) d\mu
(g)$, where $\psi$ is a bounded Borel function.  This operator
allows to express various quantities of probabilistic significance
for the left random walk $X_n^g$.  For example, if $B$ is a Borel
subset of $G$:
$$\mP \{ X_n ^g\in B \} = P_\mu ^n  1_B(g)$$ for any $g \in G$.  In
particular, if $\mu$ is adapted and recurrent then for any positive
continuous function $\phi$ on $G$  and any $g
\in G$: $\sum _0 ^\infty P_\mu ^k \phi (g) = +\infty$.  For the sake
of completeness we will give proofs for certain general properties
used below and some of which have been considered in the context of
Markov operators on measured spaces (see \cite{Fo}).  We will also
use the framework of skew products in the context of $G$-spaces (see
\cite{Fur}, \cite{Fu} and \cite{Gu3}).

A positive Borel function $f$ on $G$ is said to be left $\mu$-harmonic
(resp. left $\mu$-superharmonic) if $P_\mu f =f$ (resp. $P_\mu f \leq f$).
If $\mu$ is adapted and recurrent then continuous positive superharmonic
functions are constant (see below).

A useful concept when dealing with recurrent random walks is that of
induced random walk.  If $\mu$ is recurrent and $U$ is an open subgroup
the return time $T$ to $U$ is defined by
$T = \inf \{ n \geq 1 \mid X_n ^e \in U \}$ and we denote by $\mu ^U$ the
law of $X_T^e$.  Then $\mu ^U$ is clearly a recurrent probability on
$U$.  In particular, if $G$ is recurrent, then any open subgroup of $G$ is
also recurrent.  Also, if $H$ is a closed normal subgroup of $G$ and
$\mu$ is recurrent, then the projection of $\mu$ on $G/H$ is also
recurrent.  In particular, if $G$ is recurrent, then $G/H$ is also
recurrent.  In dealing with the structure of recurrent groups one can
reduce the study to the case of compactly generated groups.  This is because,
if $W$ is a compact symmetric neighborhood of $e\in G$, then the subgroup
generated by $W$, that is, $U = \cup _{n \geq 0} W^n$ is open, hence one
can consider the induced random walk defined by $\mu$ which is also
recurrent.  Conversely, if $G$ is a union of compactly generated subgroups and on
each of these compactly generated subgroups, symmetric random walks with
compactly supported density are recurrent, then one can use the method in
\cite{Ra2} to construct recurrent random walks on $G$.
In particular, it suffices to prove the quadratic growth
conjecture for compactly generated groups.

A basic fact proved here is that every closed subgroup of a
recurrent group is unimodular.  This extends the known fact that a
recurrent group is unimodular, strongly used in \cite{GKR} for the
early classification of recurrent locally compact groups.  It
extends also Theorem 3.26 of \cite{Wo}.

We now give a few typical examples of recurrent groups.  The euclidean
motion group $\mG _2= O(2) \ltimes \R ^2$ of the euclidean plane
was already mentioned above.  If $C$ is a
compact group, then the groups of the form $C\times H$ with $H$ a
closed subgroup of $\R ^2$ are recurrent.  These two classes of groups arise
geometrically as
similarity groups or products of similarity groups relative to a euclidean
or ultrametric norm.  It follows from \cite{Ra2} that any $p$-adic
unipotent algebraic group is also recurrent.

If $A$ is any hyperbolic automorphism of a compact abelian group $C$, then
one can form the semidirect product $G= \Z \ltimes C$ where the action of
$\Z$ on $C$ is given by $A^n $ ( $n \in \Z$).  For example, one can take
$C= \mT ^2$ and $A$ is given by the matrix $\pmatrix {2 &1 \cr 1 &1 }$.
These groups are not isomorphic to closed subgroups of linear groups but are
recurrent.

\end{subsection}

\begin{subsection}{A class of Markov operators}

Let $E$ be a locally compact second countable space.  We denote by
$C_b(E)$ the space of bounded continuous functions on $E$, by
$C_c(E) \subset C_b(E)$, the subset of compactly supported functions
and $C_b^+(E)$ (respectively, $C_c^+(E)$) the set of positive elements
in $C_b(E)$ (respectively, $C_c(E)$). Let $M^1(E)$ be the set of
Borel probability measures on $E$.  Here, by a Markov operator on
$E$, we mean a positive operator $P$ on $C_b(E)$ such that $P1 =1$.
Then $P$ defines a transition probability $P(x, \cdot )$. Clearly
$P$ acts on $M^1(E)$. We will also consider its action on some
positive Radon measures. If $\eta$ is such a measure and if for any
$\phi \in C_c^+(E)$, $\eta (P\phi )$ is finite, then $\phi \mapsto
\eta (P\phi )$ defines a Radon measure which we will denote by
$P\eta$.  In particular $\eta$ will be said to be $P$-invariant if
$P\eta$ is defined and $P\eta = \eta$.

As a special case, we will consider a $G$-space $E$ and $P$ defined by
$$P\psi (x) =\int \psi (gx) d\mu (g) $$ where $\mu \in M^1(G)$ and $\psi \in
C_b(E)$.  In this case the trajectories, starting from $x \in E$, for the
associated Markov chain can be written as $$X_n(\oa ) x= Y_n (\oa )\cdots
Y_1(\oa ) x ~~ (n >0), ~~X_0(\oa )x = x$$ with
$\oa \in \Oa$.  Then for $u \in C_b^+(E)$, we have $$\sum _0 ^\infty
P^ku(x) =
\int \sum _0 ^\infty u(X_n(\oa )x) d\mP (\oa ).$$

The first part of
the following was proved in \cite{L}.  In view of the role of this
result here and in other contexts we provide a proof different from
\cite{L}, we give examples and complements; see also \cite{Pe} for other
examples of probabilistic significance.

\bp\label{p1} Let $E$ be a locally compact second countable space and $P$ be a
Markov operator on $E$.  Assume that there exists  $u \in C_c^+(E)$ such that
for any $x \in E$, $\lim _{n \to \infty } \sum _0 ^n P^k u (x) = \infty$.  Then
there exists a $P$-invariant positive Radon measure $\nu$ on $E$ with $\nu (u)
>0$.  If $\nu$ is unique up to normalization, then we have the following
convergence: $$\lim _{n \to \infty } {\sum _0 ^n P^k \phi (x) \over \sum _0 ^n
P^k u(x)} = {\nu (\phi )\over \nu (u)}$$ for all $x \in E$ and all $\phi \in
C_c(E)$.
\ep

\bo
Let $\phi \in C_c^+(E)$.  Since the sequence $\sum _0 ^n P^k u $ is
increasing tending to $+\infty$ and the support of $\phi$ is compact,
there exists $r \in \N$
such that $\phi \leq \sum _0 ^{r} P^k u$.  It follows that for any $n
\in \N$, $\sum _0 ^n P^k \phi \leq r \sum _0 ^n P^k u +r ^2
||u||_\infty$.  In particular, if $x \in E$ is fixed, then
$$(\sum _0 ^n P^k \delta _x) (\phi ) \leq r \sum _0 ^n (P^k \delta _x)
(u) + r^2 ||u||_\infty.$$  It follows that for any $\rho \in M^1(E)$, if
$\rho _n = \sum _0 ^n P^k \rho$, then $$\rho _n (\phi ) \leq r \rho _n (u)
+r^2 ||u||_\infty.$$ Since $\lim _{n \ra \infty} \rho _n (u) =
+\infty$, for $n$ large we have $\rho _n (u) >0$, hence the sequence
${\rho _n (\phi ) \over \rho _n (u)} = \eta _n (\phi)$, say, is bounded.
This implies that the sequence of Radon measures $\eta _n$ is relatively
compact in the weak topology.  On the other hand $P\eta _n = \eta _n +
{P^n\rho -\rho \over \rho _n (u)}$.  Since for any $\phi \in C_c^+(E)$,
$0\leq P^n\rho (\phi ) \leq ||\phi ||_\infty$, the sequence $\epsi _n =
{P^n \rho - \rho \over \rho _n (u)}$ converges to zero weakly.  By weak
compactness we can extract a convergent subsequence $\eta _{n_k}$ of $\eta
_n$ such that $\lim _{k \ra \infty} \eta _{n_k} = \nu $ and $\nu (u) =1$.
Then for $\phi \in C_c^+(E)$, the relation $P\eta _n(\phi ) = \eta _n
(\phi ) +\epsi _n (\phi )$ gives
$\lim _{k \ra \infty } \eta _{n_k}(P\phi ) = \nu (\phi )$.
Let $\phi , \psi \in C_c^+(E)$ with
$\psi \leq P\phi$.  Then the relation  $\eta _n(P\phi ) = \eta _n
(\phi ) +\epsi _n (\phi )$ implies  $\eta _n(\psi )\leq  \eta _n
(\phi ) +\epsi _n (\phi )$.  Hence, in the limit $\nu (\psi )
\leq \nu (\phi )$.  Since $P\phi$ is an increasing limit of elements of
$C_c^+(E)$, we get that $P\phi$ is $\nu$-integrable and
$\nu (P\phi ) \leq \nu (\phi )$.  Since $\phi$ is arbitrary, $P\nu \leq
\nu$.  In order to show $P\nu = \nu$, we consider the positive Radon
measure $\eta = \nu - P\nu$ and we observe that $\sum _0 ^n P^k \eta
= \nu - P^{n+1} \nu$.  If $\eta \not = 0$, then the condition
$\lim  _{n \ra \infty} \sum _0 ^n P^k u = +\infty$ implies that $\lim _{n
\ra \infty} (\sum _0 ^n P^k \eta )(u) = +\infty$.  This contradicts the
fact that $\sum _0 ^n P^k \eta (u)$ is bounded by $\nu (u)$.   Hence
$P\nu = \nu$.

If the Radon measure $\nu$ is uniquely defined by $P\nu = \nu$, $\nu
(u) =1$, we can improve the above considerations: every limit point
$\eta$ of $\eta _n$ satisfies $P\eta = \eta$ with $\eta (u) =1$,
hence $\eta = \nu$ and so by compactness, the convergence to $\nu$
of $\eta _n = {\rho _n \over \rho _n (u)}$ follows. \eo

\br\label{rk1}
(a) The condition $\sum _0 ^\infty P^k u = +\infty$ on $E$ is
satisfied for some $u\in C_c^+(E)$ if there exists a relatively compact open set
$U\subset E$ such
that for any $x \in E$ the trajectories of the Markov chain defined by $P$
visit $U$ infinitely often with positive probability.  In particular, if $E=G$ and
$P_\mu = P$ is associated to an adapted recurrent $\mu \in M^1(G)$, then
the condition that $\sum _0 ^\infty P^k u = +\infty$ is satisfied.
Also in this case $\nu$ is unique up to normalization and equal to a Haar
measure on $G$.

(b) We will use the above result in the following situation.  $E$ is a
$G$-space, $\mu \in M^1(G)$ is recurrent and adapted and $P$ is defined by
$P\psi (x) = \int \psi (gx) d\mu (g)$ where $\psi \in C_b(E)$.  The
condition that $\sum _0 ^\infty P^k u = +\infty$ is satisfied if there
exists a compact subset $C$ of $E$ with $GC = E$.  However the proposition
can be used in various situations not related to recurrence of $G$, where
$E$ is a non-compact homogeneous space of $G$.  Some classes of examples
are discussed in section 5.
\er

The next result is well-known if $E$ is a countable discrete space.

\bp\label{p2}
Let $E$ be a locally compact second countable space and $P$ be a Markov
operator on $E$ which satisfies
$$\sum _0^\infty P^k u = +\infty {\rm ~~on ~~} E$$  for all $u\in C_c^+(E)$.
If $f$ is a positive continuous function on $E$ with satisfies
$Pf\leq f$, then $f$ is constant.  Any non-zero $P$-invariant measure $\nu$
satisfies $\nu (\phi ) >0$ for any $\phi \in C_c^+(E)$.
\ep

\bo
The function $\psi = f-Pf$ is a positive continuous function such that
$\sum _0 ^n P^k \psi = f- P^{n+1}f \leq f$.  If $\psi$ is not zero, then
the assumptions on $P$ implies $\sum _0 ^\infty P^k \psi = +\infty$ on $E$
which contradicts the boundedness of $\sum _0 ^n P^k \psi$.  Hence we have
$Pf = f$.

Let $r \geq 0$.  Then function $f_r = \inf (r, f)$ satisfies
also $Pf_r \leq f_r$.  The above result applied to $f_r$ gives $Pf_r =
f_r$.  But for any $r < \sup _{x\in E} f(x)$, the functions $f_r$ has
maximum $r$.  We consider the closed set $E_r = \{ x\in E \mid f_r (x) =
r \}$.  Then if $x \in E_r$, the equation $Pf_r = f_r$ implies that $P(x,
E_r) = 1$.  In other words, $E_r$ is $P$-invariant.  If $f$ takes two
distinct values $r\not =s$, we can find $u \in C_c^+(E)$ with $u = 0$ on $E_s$.
Then for any $n \in \N$, $P^n u(x) = 0$ if $x \in E_s$ and hence $\sum _0 ^\infty
P^ku(x) = 0$ on $E_s$.  This contradicts the assumption that $\sum _0^\infty P^k
u = +\infty$ on $E$.  Hence $f$ is constant on $E$.

Let $\phi \in C_c^+(E)$.  Then since
$\sum _0 ^\infty P^k \phi = +\infty$ on $E$, given any $u\in C_c^+(E)$, there
exists $m \in \N$ such that $u \leq \sum _0 ^m P^k\phi$.  Then $\nu (u) \leq
\sum _0 ^m \nu (P^k \phi ) = (m+1)\nu (\phi )$.  Now choose $u \in C_c^+(E)$
so that $\nu (u) =1$, hence we have $\nu (\phi ) \geq {1\over m+1} >0$.
\eo
\end{subsection}

\begin{subsection}{ Recurrence and ergodicity for stationary measures}

Here $E$ is a locally compact second countable $G$-space,
and $P$ is defined by $$P\psi (x) = \int
\psi (gx) d\mu (g)$$ where $\mu \in M^1(G)$.  We consider also a
positive Radon measure $\lam$ on $E$ which is $P$-invariant.  In
general $\lam$ will be of infinite mass but we will deal here with
conservativity conditions which will allow to reduce the situations
to the finite mass case.  For $(\oa , x) \in \Oa \times E$ we write
$X_n(\oa ) x = Y_n (\oa ) \cdots Y_1(\oa ) x$ and we endow $\Oa
\times E$ with the measure $\mP \otimes \lam$.

\bd Let $(E, P, \lam)$ be as above.  We say that $(E, P, \lam )$ has
property $R$ if for every open relatively compact set $U\subset E$,
$\mP \otimes \lam $-a.e., $(\oa ,x) \in \Oa \times U$, there exists
$n = n(\oa , x)\in \N$ such that $X_n(\oa )x \in U$. \ed

For example if $\mu$ is a left random walk on a locally compact group $G$
and $\lam = \lam _G$, then the triple $(G, P, \lam)$ has
property $R$ if and only if $\mu$ is recurrent.

\br
In the definition of Property $R$ we could have restricted $U$ to belong to a
family of relatively compact open sets $U_k$ ($k \in \N$) covering the support
of $\lam$, which may be seen as follows.   Replacing $U_k$ by $\cup _{i\leq k}
U_k$, we may assume that $(U_k)$ is increasing.  Let $U$ be any relatively
compact open set.  Then there exists a $k\in \N$ such that
$\lam (U\setminus U_k) =0$, hence for $\mP \otimes \lam $-a.e. $(\oa ,x) \in \Oa
\times U$, there exists $n \in \N$ such that $X_n(\oa )x \in U$.
\er

We consider the skew product $(\Oa \times E, \tilde \theta , \tilde
\lam )$ where $\tilde \lam = \mP \otimes \lam$ and the map $\tilde
\theta $ is defined by $\tilde \theta (\oa , x) = (\theta \oa , Y_1
(\oa ) x)$.  We observe that $\tilde \lam$ is $\tilde
\theta$-invariant since for any $\phi \in C_c(E)$ and a Borel
function $\psi$ on $\Oa$, we have

$$\begin{array}{lcl} \tilde \theta \tilde \lam (\psi \otimes \phi)
& = & \int \psi (\theta \oa) \phi (Y_1(\oa)x)d\mP(\oa ) d\lam (x)
\\ & = & \int \psi (\oa ') \phi (gx) d\mP (\oa ') d\mu (g) d\lam (x) \\
& = & \int \psi (\oa ') \phi (y) d\mP(\oa ') d\lam (y) \\
& = & \tilde \lam (\psi \otimes \phi). \end{array}$$ Property $R$
allows us to define a return time $\tau$ and an induced
transformation $\tilde \theta _{\Oa \times U}$ on $\Oa \times U$
where $U$ is an open relatively compact set with $\lam (U)>0$ by
$$\tilde \theta _{\Oa \times U} (\oa , x) = (\theta ^{\tau (\oa, x)} (\oa ),
X_{\tau (\oa, x)} x) = \tilde \theta ^{\tau (\oa ,x)}(\oa ,x)$$
where $\tau (\oa ,x) = {\rm Inf} \{ n \geq 1 \mid X_n(\oa )x \in U \}$.

We recall that, for a dynamical system $(X, S, \eta )$ a measurable subset $D$ of
$X$ is said to be wandering if $S^{-k}D \cap D = \emptyset$ for any $k \geq 1$.
Then using the Poincar\'e recurrence theorem for $\tilde \theta _{\Oa \times U}$
with $\lam (U )<+\infty$ and standard arguments we obtain the following:

\bp\label{pr}
Let $(E, P, \lam )$ and $(\Oa \times E, \tilde \theta , \tilde \lam )$ be
as above.  Then the following are equivalent:

\be
\item [(a)] Property $R$ is valid for $(E, P, \lam )$;

\item [(b)] For any measurable $A\subset E$ with $0 <\lam (A) < +\infty$,
$\sum _0 ^\infty 1_{\Oa \times A} (\tilde \theta ^k(\oa,  x)) =
+\infty$ $\tilde \lam$-a.e on $\Oa \times A$;

\item [(c)] $(\Oa \times E, \tilde \theta , \tilde \lam )$ has no
wandering set of positive measure.
\ee
\ep

Using Proposition \ref{pr}, we see that property $R$ for $(E, P, \lam )$ is
equivalent to conservativity of the system $(\Oa \times E, \tilde \theta ,
\tilde \lam )$ (see \cite{Fo}).

Also, property $R$ is valid in the setting of Proposition \ref{p1} as shown
below.

\bc\label{cr}
Assume that for some $u\in C_c^+(E)$, $\sum _0 ^\infty P^k u = +\infty$ on
$E$, and let $\lam$ be any $P$-invariant Radon measure.  Then $(E, P, \lam )$
has property $R$.
\ec

\bo We observe that $P$ acts on $\mL ^\infty (\lam )$ and has an
adjoint operator $P^*$ on $\mL ^1 (\lam )$ defined by $$<\psi ,
P^*\phi > = \lam (\psi P^*\phi ) = \lam (\phi P\psi ) = <\phi ,
P\psi >.$$ Since $\lam$ is $P$-invariant, $P^* 1 =1$.  By
definition, since $u \in \mL ^1 _+ (\lam )$, $P^*$ is conservative.
Then it follows from Hopf's maximal ergodic Lemma (see \cite{Fo} pp.
11) that for any $\phi \in \mL ^1_+(\lam )$, $\phi >0$ $\lam$-a.e,
$\sum _0 ^\infty P^k\phi = +\infty$, $\lam$-a.e.  As a consequence
$\sum _0 ^\infty (P^*)^k \phi = +\infty$ $\lam$-a.e.  We show that
$(\Oa \times E , \tilde \theta , \mP \otimes \lam )$ has no
wandering set of positive measure.  Assume $D\subset \Oa \times E$
is such a set.  Then $\sum _0 ^\infty 1_D \circ  \theta ^k \leq 1$.
For two non-negative measurable functions $f$ and $f'$ on $\Oa
\times E$, we write
$$<f ,f '> = (\mP\otimes \lam ) (ff').$$  Also, if $v$ is a non-negative
measurable function on $E$, we write $\tilde v(\oa ,x) = v(x)$. Then
we observe for $k \in \N \cup \{ 0\}$, that $$<\tilde v , 1_D\circ
\theta ^k> = <1_D , {(P^*)^kv}>.$$  Assume $v>0$ $\lam$-a.e,
hence from above $\sum _0 ^\infty (P^*)^kv= +\infty$ $\lam$-a.e.  In
particular, since $(\mP \otimes \lam) (D)>0$, $<\tilde v , \sum _0
^\infty 1_D \circ \theta ^k> = +\infty$.  On the other hand, since
$D$ is wandering
$$<\tilde v , \sum _0 ^\infty 1_D \circ \theta ^k> \leq
(\mP \otimes \lam ) (\tilde v)= \lam (v) <+\infty $$ if $v\in \mL ^1(\lam )$.
This gives the required contradiction.
\eo

\br In case $E=G$, it follows that the condition $\sum _0^\infty \mu ^k (U) =
+\infty$ for some (hence any) neighborhood $U$ of identity is necessary and
sufficient for recurrence of $\mu$.  \er

One also considers the product space $\hat \Oa = G^ \Z$ and the corresponding
shift $\theta$.  We denote $\Oa ^-= G^{-\N}\times G$ and for $\hat \oa \in
\hat \Oa = \Oa ^-\times \Oa$, we write $\hat \oa
= (\oa ^- , \oa )$ with $\oa ^- \in \Oa ^-$, $\oa \in \Oa$.  For $k
\in \Z$ we denote $Y_k(\hat \oa )$ the $k$-th component of $\hat \oa$.
We consider the transformation $\hat \theta$ on $\hat \Oa \times E$ defined by
$\hat \theta (\hat \oa , x) = (\theta \hat \oa , Y_1(\oa ) x)$ and
we observe that the system $(\Oa \times E, \tilde \theta)$ is a factor of the
invertible system $(\hat \Oa \times E, \hat \theta )$ relative to the map
$(\hat \oa ,x ) \mapsto (\oa , x)$.

The following result will play a basic role in the study of random walks on
$G$-spaces.

\bp\label{ei}
Let $E$ be a locally compact second countable $G$-space, $\mu \in M^1(G)$ and $\lam$
be a $\mu$-stationary measure on $E$.  Then with the above notations, there
exists a unique $\hat \theta $-invariant measure on $\hat \Oa \times E$ which
has projection $\tilde \lam = \mP\otimes \lam$ on $\Oa \times E$.
\ep

\bo For the existence result we will use weak convergence of Radon
measures.  Let $\ol G = G\cup \{ \infty \}$ be the Alexandrov
compactification of $G$ and let $\hat \Oa ^\infty$ be the compact
metric space $\ol G ^\Z$.  Then $\hat \Oa ^\infty \times E$ is 
locally compact and we have a well defined continuous projection of
$\hat \Oa ^\infty \times E$ on $\ol G ^{\N -n} \times E$ for every
$n \in \N$.

Let $\Oa ^e = \{ \hat \oa \in \hat \Oa \mid \oa _k =e ~~{\rm if} ~~
k\leq 0 \}$, hence $\theta ^n(\Oa ^e) = \{ \hat \oa \in \hat \Oa
\mid \oa _k = e ~~{\rm if}~~ k\leq -n \}$.  For each $\oa \in \Oa$,
we define $\oa ^e \in \Oa ^e$ by $\oa ^e _k  = e$ if $k \leq 0$ and
$\oa ^e_k = \oa _k$ if $k>0$.  This allows to identify $\Oa \times
E$ with $\Oa ^e \times E \subset \hat \Oa \times E$ using the map
$(\oa ,x) \mapsto (\oa ^e, x)$.  The corresponding push-forward of
$\tilde \lam$ is denoted by $\tilde \lam ^e$ and we consider the
sequence $\hat \theta ^n(\tilde \lam ^e)$.  Then $\hat \theta ^n
(\Oa ^e \times E)$ is a Borel subset of the locally compact space
$\hat \Oa ^\infty \times E$ and $\hat \theta ^n (\tilde \lam ^e)$
can be considered as a Radon measure on $\hat \Oa ^\infty \times E$
which gives full measure to $\theta ^n (\Omega ^e) \times E$.
By definition:
$$\hat \theta ^n(\tilde \lam ^e) = \int \delta _{\hat \theta ^n \oa ^e} \otimes
X_n(\oa ) \lam d\mP(\oa ) ,$$  hence the projection $\hat \theta ^n
(\tilde \lam ^e)$ on $\Oa ^e \times E$ is $\int \delta _{\oa ^e}
\otimes X_n(\oa )\lam d\mP(\oa ) = \mP ^e \otimes \lam$. Then for
any function of the form $\psi \otimes \phi$ with $\psi \in C^+(\hat
\Oa ^\infty )$, $\phi \in C_c^+(E)$, $\hat \theta ^n (\tilde \lam
^e) (\psi \otimes \phi )$ is bounded by ($\sup _{\hat \oa \in \hat
\Oa ^\infty} \psi (\hat \oa )) \lam (\phi )$.  It follows that the sequence of
Radon measures $\hat \theta ^n (\tilde \lam ^e)$ on $\hat \Oa ^\infty\times
E$ is weakly relatively compact.  If $\hat \lam$ is a limit of a
subsequence, then the projection of $\hat \lam$ in 
$\prod _{i>-n}\ol G \times E \subset \ol G ^{\Z}\times E$ 
is equal to $\hat \theta ^n (\tilde \lam ^e)$ 
for every $n \in \N \cup \{ 0 \}$.  Hence $\hat
\lam$ gives full measure to $\hat \Oa \times E$.  It follows that
$\hat \lam$ is $\hat \theta$-invariant and $\hat \theta ^n (\tilde \lam
^e)$ converges weakly to $\hat \lam$.  The uniqueness of $\hat \lam$
as in the proposition follows. \eo

\br\label{nrk} In some cases we can get a more explicit form of
$\hat \lam$.  We denote by $\hat \oa ^e$, the projection of $\hat
\oa \in \hat \Oa$ on $\hat \Oa ^e$, hence:

$$\tilde \lam ^e = \int \delta _{\hat \oa ^e} \otimes \lam d\hat \mP (\hat \oa )$$
$$\hat \theta ^n (\tilde \lam ^e) = \int \delta _{\theta ^n \hat \oa ^e} \otimes
X_n(\hat \oa )\lam d\hat \mP (\hat \oa ) = \int \delta _{\theta ^n
(\theta ^{-n}\hat \oa ^e)}\otimes Y_0(\hat \oa ) \cdots
Y_{-n+1}(\hat \oa )\lam d\hat \mP(\hat \oa ).$$  Since $\lim _{n \to
\infty } \theta ^n(\theta ^{-n}\hat \oa ^e) = \hat \oa$ in $\hat \Oa
^\infty$, it follows that for any $\psi \in C^+(\hat \Oa ^\infty )$,
the sequence of measures $\int Y_0 (\hat \oa ) \cdots Y_{-n}(\hat
\oa )\lam \psi (\hat \oa ) d\hat \mP(\hat \oa)$ converges weakly to the
Radon measure $\lam _\psi$ defined by $\lam _\psi (\phi )= \hat \lam 
(\psi \otimes \phi)$. This fact is well known in the situations 
of $\mu$-boundaries (see \cite{Fu}) where $E$ is compact 
and $\lam \in M^1(E)$.  Then $Y_0(\hat \oa ) \cdots
Y_{-n}(\hat \oa ) \lam$ converges $\hat \mP$-a.e. to $\lam _{\hat \oa } \in
M^1(E)$ and $\hat \lam = \int \delta _{\hat \oa }\otimes \lam _{\hat
\oa }d\hat \mP(\hat \oa)$. Also if $\lam$ is $G$-invariant $Y_0(\hat
\oa )\cdots Y_{-n}(\hat \oa )\lam = \lam$, hence $\hat \lam = \hat
\mP\otimes \lam$.  If $E$ is a $\mu$-boundary $\lam _{\hat \oa }$ is
a Dirac measure and hence $\hat \lam \not = \hat \mP \otimes \lam$. Another
situation where $\hat \lam$ can be calculated will be considered in
Section 5, example 3. \er

\bp\label{ec}
With the above notations, assume that $(E, P, \lam )$ satisfies the
conditions:

\be
\item [(a)] property $R$ is valid;

\item [(b)] the condition $Pf=f$, $f \in \mL ^\infty (\lam )$ implies $f$
is constant.
\ee
Then the system $(\hat \Oa \times E, \hat \theta , \hat \lam)$ is
ergodic.  The converse is valid, except if $E= \Z$ and $\mu = \delta _g$ acts
by translations on $\Z$.
\ep

The proof depends on the following lemma:

\bl\label{fl}
Let $\cH$ be a Hilbert space and $\cL$ be a closed subspace and $S$ be a
contraction of $\cH$ such that $$S(\cL) \subset \cL, ~~ \overline {\cup
_{n\geq 0} S^{-n} (\cL) }= \cH .$$  Assume that the restriction of $S$ to
$\cL$ has a unique fixed point $\phi \in \cL$.  Then $\phi$ is the unique
fixed point of $S$ in $\cH$.
\el

\bo Assume $f \in \cH$ satisfies $Sf =f$.  Let $f_n$ be the
orthogonal projection of $f$ onto $\cH _n = S^{-n} (\cL)$.  Since
$Sf_n \in \cH _n$ and $Sf = f$, we get that $||f-f_n || \leq
||f-Sf_n|| = ||S(f-f_n)|| \leq ||f-f_n||$ as $S$ is a contraction.
Thus, $||f-f_n|| = ||f-Sf_n||$.  This implies from the definition of
the orthogonal projection that $f_n=Sf_n$. Now, $f_n = S^nf_n \in
\cL$.  It follows from the uniqueness of $\phi$ that $f_n = \phi$.
Since $\ol {\cup _{n\ge 0} S^{-n}(\cL) }= \cH$, we have
$\phi = \lim _{n \ra \infty} f_n = f $. \eo

For any Polish space $X$ we denote by $\cB (X)$ the $\sigma$-algebra of its
Borel subsets.

\bo $\!\!\!\!\!$ {\bf of Proposition \ref{ec}}\ \ To begin with we
show that condition (b) implies ergodicity of the non-invertible
system $(\Oa \times E, \tilde \theta , \tilde \lam )$.  Let $\cA
\subset \cB(\Oa \times E)$ be the $\sigma$-algebra of sets of form
$\Oa \times B$ with $B\in \cB(E)$ and $\cA _n \subset 
\cB(\Oa\times E)$ be
the $\sigma$-algebra defined by the coordinates $Y_k$ $(1\leq k \leq
n)$.  Assume $F \in \mL ^\infty (\mP \otimes \lam )$ is $\tilde
\theta$-invariant and denote $f = \mE(F|\cA)$. Then $f(x) = \int
F(\oa, x) d\mP (\oa ) \in \mL ^\infty (\lam)$ and $Pf =f $. Also,
$\mE (f|\cA \vee \cA _n) (\oa , x) = f(X_n(\oa ) x)$ and the
martingale convergence theorem implies that $F(\oa , x) = \lim _{n
\to \infty} f(X_n(\oa )x)$, $\tilde \lam$-a.e.  Condition (b)
implies that $f$ is constant, hence $F$ is also constant.

Let $U \subset E$ be open relatively compact set, $A=\Oa \times U$
and $\hat A = \Oa ^-\times A \subset \hat \Oa \times E$.  Condition
(a) implies that for $(\Oa \times E, \tilde \theta , \tilde \lam)$,
the return time $\tau (\oa ,x)$ from $A$ to $A$ is defined a.e on
$A$. Let $\tilde \theta _A$ be the induced transformation on $A$.
Then the induced transformation on $\hat A$ is well defined by $\hat
\theta _{\hat A} (\hat \oa , x) = \hat \theta ^{\tau (\oa , x)}
(\hat \oa , x)$. The restriction $m_A$ (resp. $m_{\hat A}$) of $\tilde \lam$ 
(resp. $\hat \lam)$ to $A$ (resp. $\hat A$) is $\tilde \theta _A$-invariant 
(resp. $\hat \theta _{\hat A}$-invariant) and $(A, \theta _A , m _A)$ is a 
factor of $(\hat A, \hat\theta _{\hat A} , m_{\hat A})$.  We show ergodicity of 
the corresponding systems as
follows.  Let $C\in \cB (A)$ with $\tilde \theta _A^{-1} (C) =C$,
$\mP \otimes \lam (C) >0$, and $C' = \cup _{k\geq 0} \tilde \theta
^{-k} (C)$.  Since $\tilde \theta ^{-1} (C') \subset C'$ and $(\Oa
\times E , \tilde \theta , \tilde \lam )$ has no wandering set of
positive measure, one has $\tilde \theta ^{-1} (C') = C'$ mod
$\tilde \lam$, hence $C' = \Oa \times E$ by ergodicity of $(\Oa
\times E, \tilde \theta , \tilde \lam )$.  Then $C= C ' \cap A = A $
mod $\tilde \lam$. Now Lemma \ref{fl} gives the ergodicity of $\hat
\theta _{\hat A}$ which may be seen as follows. Since $m_A$ is the
projection of $m_{\hat A}$ on $\hat A$ one can take $\cH = \mL ^2
(m_{\hat A})$, $\cL = \mL^2 (m_A)$, $S= \hat \theta _{\hat A}$.
Since $\vee _{-\infty} ^{+\infty} \hat \theta ^k (\cA) = \cB (\hat
\Oa \times E)$, one has $\cup _0^{+\infty} \hat \theta _{\hat
A}^{-k} (\cA \cap \hat A) =\cB (\hat A)$, hence $\cup _0 ^\infty
S^{-k} (\cL) = \cH$.  Also the restriction of $S$ to $\mL ^2(m_A)$
is $\tilde \theta _A$.  Since $\tilde \theta _A$ is $m_A$-ergodic,
Lemma \ref{fl} implies $S= \hat \theta _{\hat A}$ is also $m_{\hat
A}$-ergodic. Finally the ergodicity of $\hat \theta$ is obtained as
follows.  Let $\hat C \in \cB(\hat \Oa \times E)$ with $\hat \theta
^{-1}(\hat C) = \hat C$, $(\hat \lam ) (\hat C) >0$. Then $\hat C
\cap \hat A$ is $\hat \theta _{\hat A}$-invariant, hence by
ergodicity $\hat C \cap \hat A = \hat A$, $\hat C \supset \hat A$.
Since $E$ is a union of relatively compact open sets $U_n$ with $0<
\lam (U _n) $, one gets $\hat C = \hat \Oa \times E$ mod $\hat \mP
\otimes \lam$.

For the converse, we observe that if $\mu$ is not a Dirac measure,
then $(\hat \Oa , \hat \lam)$ is non-atomic.  Furthermore, if $\mu =
\delta _g$, ergodicity of $\hat \theta $ on $\hat \Oa \times E$ is
equivalent to ergodicity of the action of $g$ on $E$. Also, in this
case $(E, \lam )$ atomic is equivalent to $(\hat \Oa \times E, \hat
\lam )$ atomic. Furthermore, $(E, \lam )$ atomic and $g$ ergodic
implies that $(E, g, \lam )$ reduces to the translation on $\Z$.  In
the opposite case, ergodicity of $\hat \theta$ implies that $(\hat
\Oa \times E, \hat \theta , \hat \lam )$ has no wandering set.  Then
by Proposition \ref{pr} property $R$ is valid for $(E, P, \lam )$.
If $f \in \mL ^\infty (\lam )$ satisfies $Pf=f$, one takes a Borel
version of $f$ again denoted by $f$, which satisfies $Pf =f$.  Then
for every $x \in E$, the sequence $(f(X_n(w)x))$ is a bounded
martingale with respect to $\mP$ and the natural filtration on
$\Oa$.  By the martingale convergence theorem we have
$$F(\oa, x) = \lim _{n\to \infty} f(X_n(\oa )x)$$ $\mP$-a.e and $f(x) = \int F(\oa ,x)
d\mP (\oa )$.  Also $F(\oa , x)$ is $\tilde \theta $-invariant mod $\mP
\otimes \lam$.  Then the ergodicity of
$(\hat \Oa \times E , \hat \theta , \hat \lam )$
implies that $F$ is constant $\mP \otimes \lam$-a.e, hence $f$ is constant
$\lam$-a.e.
\eo

Using the above remark \ref{nrk} and Proposition \ref{p2} we get the
following well-known result.

\bc\label{co2}
Assume that if $\mu \in M^1(G)$ is adapted and recurrent.  Then the
skew-product $(\hat \Oa \times G, \hat \theta , \hat \mP \otimes \lam _G)$ 
is ergodic.
\ec

We now relate ergodicity of $(\hat \Oa \times E , \hat \theta , \hat \lam )$
to extremality of $\lam$.

\bc\label{co3} Assume $\lam$ is a $P$-invariant Radon measure on $E$
and $(E, P, \lam)$ satisfies property $R$.  Then $(\hat \Oa \times E
, \hat \theta , \hat \lam )$ is ergodic if and only if
$\lam$ is extremal. \ec

\bo Assume $\lam$ is extremal and let $U$ be a relatively compact
open subset of $E$ with $\lam (U) >0$.  In order to prove ergodicity
of $\hat \theta$, we consider the return time of $U$ defined by
$$H(\oa ,x) = \inf \{ n \in \N \mid X_n(\oa )x \in U \}, ~~~~ (\oa ,
x) \in \Oa \times U .$$ Also, we can define the induced Markov
operator $P^H$ of $P$ on $U$.  Let $f$ be a bounded Borel function
such that $Pf=f$. Since $H(\oa ,x)$ is a stopping time we have
$f(x) = \mE (f(X_{H(x,\oa )}(\oa ) x))= P^Hf(x)$.  Let $\lam _U =
1_U\lam$ be the restriction of $\lam$ to $U$, hence $P^H\lam _U =
\lam _U$.  Since $\lam _U$ is finite and $f1_U$ is $P^H$-invariant,
we have $P^H(f\lam _U)=f\lam _U$.  Then by the classical
regeneration method we can construct a $P$-invariant measure $\rho
_f$ on $E$ such that the restriction of $\rho _f$ to $U$ is $f\lam
_U$.  We recall the construction: let $\rho (\oa , x) $ be the random
measure given by $\rho (\oa ,x ) = \sum _0 ^{H(\oa , x)-1} \delta
_{X_k(\oa )x}$ and observe that if $$Y_1(\oa ) \rho (\oa , x ) = \rho
(\oa, x) +\delta _{X_{H(\oa, x)}x} -\delta _x$$ and $\rho _f$ is
defined by $\rho _f = \mE (\int \rho (\oa, x ) f(x) d\lam _U(x))$,
its $P$-invariance follows from the $P^H$-invariance of $f\lam _U$.
If $f =1$, one may verify that this procedure gives $\rho _f
=1_{U'}\lam$ with $U'= T_\mu U$ where $T_\mu$ is the closed
semigroup generated by the support of $\mu$.  Since $|f|$ is bounded
by $c>0$, it follows that $\rho _f$ is a Radon measure with $\rho _f
\leq c\lam$ and $P\rho _f = \rho _f$.  The extremality property of
$\lam$ gives that $\rho _f$ is proportional to $\lam$, in particular
$f\lam _U$ is proportional to $\lam _U$, that is, $f1_U $ is
constant $\lam _U$-a.e.  Since $U$ was arbitrary with $\lam (U) >0$,
we can conclude $f$ is constant $\lam$-a.e.  Now applying
Proposition \ref{ec} to $(E,P, \lam)$, we get the ergodicity of $(\hat \Oa
\times E , \hat \theta , \hat \lam )$.

Conversely, let $\lam '$ be another positive Radon measure with
$P\lam ' = \lam '$, $\lam ' = f\lam$ where $f \in \mL ^\infty (\lam
)$.  We will now claim that $f$ is constant.  As above we write
$\lam _U = 1_U \lam $, $\lam _U ' = 1_U\lam '$.  We observe that
since property $R$ is valid $H(\oa ,x)$ is well defined for $\mP
\otimes \lam $-a.e $(\oa , x) \in \Oa \times U$.  Then $P^H \lam _U
' = \lam _U'$, $P^H \lam _U = \lam _U$ and $\lam _U' = (f1_U) \lam
_U$.  It follows that $\mP \otimes \lam _U$ is $\tilde \theta
^H$-invariant.  By ergodicity of $(\Oa \times U, \tilde \theta ^H ,
\mP\otimes \lam _U )$: $f1_U $ is constant $\lam _U$-a.e. Since $U$
is arbitrary with $\lam (U)
>0$, we get $f$ is constant $\lam $-a.e and $\lam '$ is proportional
to $\lam$. \eo

\end{subsection}

\begin{subsection}{Recurrent groups and contraction subgroups}

Let $G$ be a locally compact group.  For $g \in G$, we define the
contraction subgroup $C_g$ of $g$ by $$C_g = \{ x \in G \mid
\lim _{n \ra \infty } g^n x g^{-n} = e \}.$$  Then $C_g$ is a subgroup of
$G$ normalized by $g$.  In general, $C_g$ is not closed.  However, if
$G = GL(n , \mF)$ for some local field $\mF$, then $C_g$ is an algebraic
$\mF$-subgroup of $G$, hence $C_g$ is closed.  If $G$ is a closed subgroup
of $GL(n, \mF)$, then $C_g$ is closed in $G$.  As a simple consequence we obtain
that if $G$ is a
closed subgroup of $\prod _{i \in I} GL(d_i, \mF_i)$ (finite $I$), $C_g$
is closed.  We now prove the following basic lemma on $C_g$.

\bl\label{l1}
Let $G$ be a locally compact group.  Assume that $g\in G$ is such that
$C_g$ is closed and $C_g \not = (e)$.  Then the closed subgroup $L$
generated by $g$ and $C_g$ is not unimodular.
\el

\bo
Let $H= C_g$.  Then $H$ is a normal subgroup of $L$ and hence $L/H$ is
abelian.  For a locally compact group $M$, we recall that $\lam _M$ is a left
Haar measure on $M$.  Let $\phi \in C_c(L)$.  We observe that the formula
$$\lam _L(\phi ) =\int d\lam _{L/H}(\overline {u}) \int \phi (uh)d\lam
_H(h)$$ defines a left Haar measure on $L$.  For any $x\in L$,  the map $h
\mapsto xhx^{-1}$ is an automorphism of $H$ and we denote by $\Delta (x)$
its module.  Furthermore, for $\phi \in C_c(L)$,  we define
$\phi ^x \in C_c(L)$ by $\phi ^x(u) = \phi (ux)$
for all $x, u \in L$.  Since $L/H$ is unimodular,
we see that $\lam _L (\phi ^x) = \Delta (x^{-1}) \lam _L(\phi )$
for all $x \in L$ and $\phi \in C_c(L)$.
Since $C_g$ is closed, for any compact neighborhood $W$ of $e$ in $G$,
$g^n W g^{-n} \ra e$ (cf. \cite{Wa}).  This implies that $\lim _{n \ra
\infty} \lam _H (g^n Wg^{-n}) = 0$, hence $\Delta (g) <1$.  It follows
that $\lam _L (\phi ^g) > \lam _L(\phi )$ for any $\phi \in C_c(L)$, hence
$L$ is not unimodular.
\eo

\bp\label{p3}
Let $G$ be a locally compact group and $E$ be a locally compact second
countable $G$-space.  Suppose there exists a compact subset
$C$ of $E$ such that $E = GC$ and $G$ is recurrent.  Then there exists a
$G$-invariant Radon measure on $E$.
\ep

\bo
Let $\mu \in M^1(G)$ be recurrent.  Let $P$ be the Markov operator on
$G$ defined by $P\phi (x) = \int \phi (gx) d\mu (g)$ for all $x \in E$.
Let $u \in C_c^+(E)$ be such that $u >0$ on $C$.   Now for any $x \in E$,
there exists $h \in G$ such that $hx\in C$, hence there exists $\delta >0$
and an open neighborhood $V_x$ of $h$ in $G$ such that
$u(gx) >\delta >0$ for all $g \in V_x$.  Then
$P^ku (x) = \int u (gx) d\mu ^k(g) \geq \delta \mu ^k (V_x)$ for all $k
\geq 1$.  Since $\mu$ is recurrent, $\sum _0 ^\infty P^ku (x) = +\infty$.
Using Proposition \ref{p1}, we get the existence of a $P$-invariant Radon
measure $\nu$ on $E$.  For any $\phi \in C_c^+(E)$, we define $h _\phi $ on
$G$ by
$h _\phi (g) = g\nu (\phi ) = \int \phi (gx) d\nu (x)$.  Since $\nu$ is
$P$-invariant, $h _\phi$ is right $\mu$-harmonic, that is,
$\int h_\phi (g g') d\mu (g') = \int \phi (gg'x) d\mu (x) d\mu (g') = \int
\phi (gy) d\nu (y) = h _\phi (g)$.  Since $h _\phi$ is continuous, positive
right $\mu$-harmonic and $\mu$ is recurrent, using Proposition \ref{p2}
for $E=G$, we get that $h _\phi$ is constant.  It follows that $g\nu
(\phi) = \nu (\phi )$ for any $g \in G$ and any $\phi \in C_c^+(G)$.
Hence $\nu$ is $G$-invariant.
\eo

\bc\label{c1}
Suppose $G$ is a locally compact recurrent group.  Then any closed
subgroup of $G$ is unimodular.  In particular, if $g \in G$ and $C_g$ is
closed, then $C_g = \{ e \}$.
\ec

\bo We first recall that for a closed subgroup $H$ of $G$, the
quotient space $G/H$ carries a $G$-invariant measure if and only if
the restriction to $H$ of the modular function of $G$ is equal to
the modular function of $H$.  We now claim that $G$ is unimodular.
For any closed subgroup $H$, let $\Da _H$ be the modular function of
$H$.  Let $g \in G$ and $H(g)$ be the closed subgroup generated by
$g$ in $G$.  Then by Proposition \ref{p3}, $G/H(g)$ carries a
$G$-invariant measure.  This implies that $\Da _G (g) = \Da
_{H(g)}(g) =1$ as $H(g)$ is abelian.  Thus, $G$ is unimodular.  Now,
let $H$ be any closed subgroup of $G$.  Then by Proposition
\ref{p3}, $G/H$ carries a $G$-invariant measure and hence $\Da _H(h)
= \Da _G(h) =1$ for any $h \in H$, as $G$ is unimodular.

Suppose $g \in G$ is such that $C_g$ is closed.  Let $L$ be the closed
subgroup generated by $C_g$ and $g$.  Then from above $L$ is unimodular and hence 
by Lemma \ref{l1}, $C_g= \{ e \}$.
\eo

In our effort to understand the structure of totally disconnected
recurrent groups, we obtain the following.

\bc If $G$ is a locally compact totally disconnected recurrent
group.  Then $G$ is uniscalar, that is, for each $g \in G$, there is
a compact open subgroup $K_g$ of $G$ such that $gK_gg^{-1}= K_g$.
\ec

\bo It follows from Corollary \ref{c1}, that closed subgroups of $G$
are unimodular and hence the result follows from Proposition 3.21 of
\cite{BW}. \eo

\end{subsection}

\end{section}

\begin{section}{Proof of the theorem \ref{mt}}

\begin{subsection}{Recurrent real Lie groups}

\bt\label{rl}
Let $G$ be a locally compact group such that $G$ has a continuous
injection into a real Lie group.  Suppose $G$ is compactly generated
and recurrent.  Then $G$ has a compact normal subgroup $K$ and a finite
index subgroup $H$ such that $K\subset H$ and $H/K$ is isomorphic to a
closed subgroup of the group $\mG _2$ of euclidean motions of the plane.
\et

\bo
We first claim that $G$ is a real Lie group.  Let $G'$ be a real Lie
group and $\phi \colon G \ra G'$ be a continuous injection.  Since $G'$ is a real
Lie group, by \cite{MZ}, there is a neighborhood $U$ of identity in $G'$ such
that $\{ e \}$ is the only subgroup contained in $U$.  Let $V= \phi ^{-1}(U)$.
Then $V$ is a neighborhood of identity in $G$.  If $N$ is a subgroup contained in
$V$, then $\phi (N)$ is a subgroup contained in $U$, hence $\phi (N)$ is trivial.
Since $\phi$ is injective, $N$ is trivial.  Now it follows from \cite{MZ} that
$G$ is a real Lie group.

Now, $G$ is a real Lie group implies that its connected component
$G_0$ is open, hence $G/G_0$ is a finitely generated recurrent
group. Thus, by \cite{VSC}, $G/G_0$ has a normal subgroup of finite
index isomorphic to ${\mathbb Z}^k$ for $k \leq 2$.  Thus, we may
assume that $G/G_0 \simeq \Z ^k$ Since $G_0$ is open, $G_0$ is a
connected recurrent Lie group. By Theorem 0.1 of \cite{Ba}, $G_0$
contains a compact subgroup $K$ and a normal vector subgroup
${\mathbb R}^l$ for $l \leq 2$ such that $G_0 \simeq K\ltimes
{\mathbb R}^l$.

If $l=0$, we get that $G_0 $ is compact and up to finite index
$G/G_0$ is a subgroup of $\Z ^2$.  We consider the case when $l=1$.
Since $K$ is connected and the only compact connected group of
automorphisms of ${\mathbb R}$ is trivial, $G_0\simeq K\times
{\mathbb R}$.  Since $G_0$ is normal in $G$, $K$ is also normal in
$G$.  Replacing $G$ by $G/K$, if necessary, we may assume that $G_0
\simeq {\mathbb R}$.  Let $\phi \colon G \ra {\rm Aut}~(G_0) \simeq
{\mathbb R}^*$ be $$\phi(g)(v) = gvg^{-1}$$ for $g \in G$ and $v \in
G_0$.  By Corollary \ref{c1}, $G$ is unimodular and since $G_0$ is
open, Haar measure on $G_0$ is $\phi (g)$-invariant for any $g \in
G$.  This implies that $|\phi (g) |=1$. Thus, $G$ contains a normal
subgroup $G_1$ of finite index such that $G_0$ is contained in the
center of $G_1$ and $G_1/G_0 \simeq {\mathbb Z}^k$.  Thus, $G_1$ is
a compactly generated nilpotent recurrent group.  By Corollary 13 of
Chapter III of \cite{GKR}, $G_1\simeq \R \times \Z$.

We next consider the case $l=2$.  The only connected compact groups
of automorphisms of ${\mathbb R}^2$ are the trivial group and
$SO_2({\mathbb R})$.  Hence, since K is connected, its conjugation
action on ${\mathbb R}^2$ factors through one of these groups.
Therefore, modulo a compact normal subgroup, $G_0$ may be assumed to
be either ${\mathbb R}^2$ or $\mG _2$, the motion group of the plane.
If $G_0 =\R ^2$, then for any $g \in G$, $g$ acts linearly on $\R
^2$ by $\tilde g$ with eigenvalues of modulus one which may be seen
as follows.  If $g \not \in G_0$, $<g>$ is isomorphic to $\Z$ and
$<g > \ltimes \R ^2$ is closed, hence unimodular by Corollary \ref{c1}.
It follows that
$|\det \tilde g| =1$.  If $\tilde g$ has an invariant line in $\R
^2$, the same argument gives that the corresponding eigenvalue has
modulus one.  Since $|\det \tilde g| =1$, the same is true of the
second eigenvalue.  If $\tilde g$ has no real eigenvalue, the
condition $|\det \tilde g| =1$ implies that $\tilde g$ has conjugate
pair of eigenvalues of absolute value one.  Now there are two
situations depending on the fact that $G/G_0$ acts on $\R ^2$
irreducibly or not.  In the first case the fact that elements
$\tilde g$ have eigenvalues of modulus one, implies that $G/G_0$
acts as a subgroup of $O(2)$.  Since $G/G_0\simeq \Z^k$, $G$ is a
two-step solvable group with relatively compact action on $[G,
G]\subset G_0$ hence by Theorem 11, Chapter IV of \cite{GKR}, we get
that $G= \R^2$.  In the second situation, $G$ has a finite index and
hence recurrent nilpotent subgroup $G_1 \supset \R ^2$.  Hence as
above, $G=G_1 = \R ^2$.   If $G_0= SO(2)\ltimes \R^2$, since the
square of every automorphism of $G_0$ is interior, we can assume by
passing to a finite index subgroup that for any $g\in G$, the inner
automorphism of $g$ restricted to $G^0$ is uniquely defined by
$g'\in G_0$.  Then the homomorphism $g\mapsto g'$ is a projection of
$G$ onto $G_0$ with kernel $N\simeq G/G_0$. Hence $G =G_0\times N$
with $N\simeq \Z ^i$, $0\leq i\leq 2$.  Then since $G$ is two-step
solvable and recurrent and $G$ acts on $\R ^2$ as a compact
subgroup, Theorem 11, Chapter IV of \cite{GKR}, gives $N$ is trivial. \eo

\bc\label{pg}
Assume $G$ is a compactly generated group of polynomial growth and is
recurrent.  Then $G$ has at most quadratic growth and has a compact normal
subgroup $K$ and a finite index subgroup $H$ such that $K\subset H$ and
$H/K$ is isomorphic to a closed subgroup of the motion group of the plane.
\ec

The corollary is a simple consequence of Theorem \ref{rl} and of the fact
that a compactly generated group $G$ of polynomial growth contains a compact
normal subgroup $K$ such that $G/K$ is a real Lie group (cf. \cite{Lo}).

\end{subsection}

\begin{subsection}{Recurrent linear groups}

We now consider closed subgroups of linear groups over local fields and
prove the quadratic growth conjecture for such groups.
We also prove Theorem \ref{mt}.

\bl\label{l2}
Let $V$ be a vector space over a non-archimedean local field $\mathbb F$
and $G$ be a subgroup of $GL(V)$.  Suppose for any
$g \in G$ and $v \in V$,  $(g ^n(v))_{n \in {\mathbb Z}}$ is
relatively compact.  Then $G$ is contained in a compact extension of a
unipotent subgroup of $GL(V)$.  In addition if $G$ is compactly generated, then
$G$ is relatively compact.
\el

\br This result is also proved in \cite{Pa}, we provide a proof as
our proof is simpler. \er

\bo
Let $V_1$ be an $G$-irreducible subspace of $V$.  Then using Burnside
density Theorem as in \cite{CG}, we conclude that the restriction of $G$ to
$V_1$ is relatively compact in $GL(V_1)$.  Using a Jordan-H\"older
sequence we get that there is a compact subgroup $K$ and an
unipotent subgroup $U$ normalized by $K$ such that
$G\subset K\ltimes U\subset GL(V)$.  Let $A$ be a compact generating set
in $G$.  Then there is a compact set $C\subset U$ with $kCk^{-1}= C$ for all
$k \in K$ (this is possible as $K$ is a compact group normalizing $U$)
such that $A\subset KC$.
Since a compactly generated subgroup of any unipotent group is relatively
compact (as the base field is non-archimedean) and $C$ is
$K$-invariant, the closed subgroup $L \subset U$ generated by $C$ is compact and
is normalized by $K$.  Thus, $G\subset K\ltimes L$ which is compact and
hence $G$ is relatively compact.
\eo

\bp\label{l3}
Let $G$ be a closed compactly generated subgroup of a linear group $GL(n, \mF)$
over a non-archimedean local field $\mathbb F$.  Suppose $C_g = \{ e \}$ for any
$g \in G$.  Then $G$ has a basis of compact open normal subgroups.
\ep

\bo Let $\Phi \colon G \ra M_n({\mathbb F})$ be given by $\Phi (x) =
x-I$ for all $x \in G$.  Then $\Phi$ is a homeomorphism onto $\Phi
(G)$ endowed with the topology induced from $M_n (\mF )$ and $\Phi
(gxg^{-1})= g\Phi (x) g^{-1}$ for all $x , g \in G$.  Let $V$ be the
smallest subspace of $M_n({\mathbb F})$ such that $V\cap \Phi (G)$
is a neighborhood of $0$ in $\Phi (G)$.  We now claim that for any
$g \in G$, $V$ is $g$-invariant and $(g^nv g^{-n})_{n \in {\mathbb
Z}}$ is relatively compact for all $v \in V$.  Since
$\mathbb F$ is non-archimedean, $G$ is a totally disconnected
locally compact group.  Then by Proposition 2.1 of \cite{JR}, there
is a basis of compact open subgroups $(K_i)$ at $e$ in $G$ such that
$gK_i g^{-1} = K_i$ for all $i \geq 1$.  Since $V\cap \Phi (G)$ is a
neighborhood of $0$ in $\Phi (G)$, there exists a $K_i$ such that
$\Phi (K_i) \subset V$.  Let $W$ be the subspace of $V$ spanned by
$\Phi (K_i)$.  Then $\Phi (K_i) \subset W\cap \Phi (G)$ is a
neighborhood of $0$ in $\Phi (G)$.  Since $V$ is the smallest such subspace
$V= W$. For any $v \in \Phi (K_i)$, $g^n vg^{-n} \in \Phi (K_i)$
for any $n \in {\mathbb Z}$ and hence $(g^nv g^{-n})_{n \in {\mathbb Z}}$
is relatively compact as $\Phi (K_i)$ is compact in
$V$.  Since $V=W$ is spanned by $\Phi (K_i)$, we get that $(g^nv
g^{-n})_{n \in {\mathbb Z}}$ is relatively compact for all $v \in
V$.  Since $g\Phi (K_i)g^{-1} = \Phi (gK_ig^{-1}) = \Phi (K_i)$ and
$W=V$ is spanned by $\Phi (K_i)$, $gVg^{-1}=V$.  Thus, $V$ is
$G$-invariant and $(g^nv g^{-n})_{n \in {\mathbb Z}}$ is relatively
compact for any $g \in G$ and $v \in V$.

Now, by Lemma \ref{l2}, $V$ contains a basis of open neighborhoods at $0$
invariant under conjugation by elements of $G$ and so $V\cap \Phi (G)$
has small $G$-invariant neighborhoods of $0$ in $V\cap \Phi (G)$.
Since $V\cap \Phi (G)$
is a neighborhood of $0$ in $\Phi (G)$, we get that $G$ contains a basis
of open invariant neighborhoods at $e$.  Since $G$ is totally
disconnected, $G$ contains a basis of compact open normal subgroups.
\eo

\bc\label{c3}
Let $G$ be a compactly generated closed subgroup of a linear group $GL(n, \mF)$
over a non-archimedean local field $\mathbb F$.  Suppose every closed subgroup of
$G$ is unimodular. Then $G$ has a basis $(K_n)$ of compact open normal subgroups.
In particular, if $G$ has polynomial growth, then furthermore, $G/K_n$ has a
finite index subgroup which is finitely generated and nilpotent.
\ec

\bo
Since $G$ is a closed subgroup of $GL(n, \mF)$, for any $g \in G$, $C_g$ is a
closed subgroup of a unipotent algebraic group, hence $C_g$ is closed.  Then
by Lemma \ref{l1}, $C_g = \{ e \}$ and first part of the corollary follows
from Proposition \ref{l3}.  If $G$ has polynomial growth, then it closed 
subgroups are unimodular (see \cite{Gu}).  So, as above we obtain using 
Lemma \ref{l1} that there is a basis $(K_n)$ of compact open 
normal subgroups.  Since $K_n$ is open, $G/K_n$ is a finitely generated group of
polynomial growth and hence the result follows from \cite{Gr}.
\eo

We now prove the quadratic growth conjecture for closed subgroups of linear
groups over local fields.

\bt\label{t1}
Let $G$ be a closed subgroup of a linear group $GL (n, \mF )$ over a local field
$\mF$.  Then $G$ is a recurrent group if and only if $G$ has at most quadratic
growth.  More precisely, if $G$ is compactly generated and recurrent, then we
have the following two cases
\be

\item[a)] when $\mF$ is archimedean, up to finite index and a modulo a
compact normal subgroup $G$ is isomorphic to a closed subgroup of the
euclidean motion group of the plane.

\item [b)] when $\mF$ is non-archimedean, up to finite index and modulo a
compact open normal subgroup $G$ is isomorphic to a subgroup of $\Z ^2 $.
\ee
\et

\bo
The implication that $G$ has at most quadratic growth implies $G$ is a
recurrent group can be proved as in Proposition 3.1 of \cite{Ra2} as explained in
section 2.1.  Here we prove the implication that $G$ is a
recurrent group implies $G$ has at most quadratic growth.

Suppose $G$ is recurrent.  In order to prove $G$ has
at most quadratic growth we may assume that $G$ is a compactly generated
group.  Now we prove the stronger statement in the second part of the
Theorem.  Since $G$ is a linear group, $C_g$ is closed in $G$ for any
$g\in G$.  Then by Corollary \ref{c1}, $C_g$ is trivial.

\noindent{\bf Non-archimedean case:}  Assume that $G$ is a linear
group over a non-archimedean field.  Then by Proposition \ref{l3}, $G$ contains
a compact open normal subgroup, say $K$.  Now $G/K$ is a finitely
generated recurrent group and hence up to finite index $G/K$ is isomorphic
to a  subgroup of $\Z ^2$ (cf. \cite{VSC}).  Since $K$ is compact, $G$
itself has at most quadratic growth.

\noindent{\bf Archimedean Case:} Suppose $\mF$ is archimedean.  Then
$\mF$ is $\R$ or $\C$.  So, we may assume that $G$ is a closed subgroup of
a linear group over $\R$.  Then it follows from Cartan's Theorem that $G$ is a
real Lie group.  Now the result is a direct consequence of Theorem \ref{rl}.
\eo

\bo $\!\!\!\!\!$ {\bf of Theorem \ref{mt}}\ \
Suppose $G$ is a compactly generated recurrent group.  For $i \in I$, let $G_i$ be
the closure of the projection of  $G$ into $GL(n_i , \mF_i)$.  Then each $G_i$
is a closed subgroup of $GL(n_i , \mF_i)$ that is compactly generated and
recurrent.

Let $J$ be the set of indices $i$ such that $\mF _i$ is non-archimedean.  Then
for $i \in J$, it follows from Theorem \ref{t1} that $G_i$ has a
compact open normal subgroup, say $K_i$.  Now $\prod _{i\in I} G_i$ is a closed
subgroup of $\prod _{i\in I} GL(n_i , \mF _i)$ and hence $G$ is a closed
subgroup of $\prod _{i\in I} G_i$.  Let $K = \prod _{i\in I} K_i$ with $K_i =(e)$
for $i \not \in J$.  Then $K$ is a compact
normal subgroup of $\prod _{i\in I} G_i$ and
$(\prod _{i\in I} G_i) /K$ is a real Lie group.  Since $K$ is compact and $G$ is
closed in $\prod _{i\in I} G_i$, $G/G\cap K$ is a closed subgroup of
$(\prod _{i\in I} G_i) /K$.  By Cartan's Theorem $G/G\cap K$ is a real Lie group.
Since $G/G\cap K$ is also recurrent, the result follows from Theorem \ref{rl}.
\eo

\br
More generally we can ask if $G$ is a compactly generated locally compact group
with a continuous embedding in $\prod _{i \in I} GL(d_i, \mF _i)$ and is
recurrent, then $G$ has at most quadratic growth.  This property is valid if all
$\mF _i$ ($i\in I$) are archimedean.  Answer to this hinges on the following:
if the image of $G$ is contained in a compact linear group over a
non-archimedean local field, is it true that $G$ is a projective limit of real
Lie groups.
\er

\end{subsection}

\end{section}

\begin{section}{Harris-recurrence}

We now consider the notion of recurrence known as Harris recurrence
(abbr. as H-recurrence) widely use in the context of Markov operators on
measured spaces (see \cite {Re}).  
We will say that $\mu$ is $H$-recurrent if for any Borel subset $B$
of $G$ of positive Haar measure, $\mP$-a.e $X_n(\oa ) \in B$
infinitely often.  If $G$ supports a $H$-recurrent $\mu$, we say
that $G$ is $H$-recurrent. It is known that $\mu$ is H-recurrent if
and only if some power of $\mu$ is non-singular and recurrent
(Theorem 4.11, Chapter 3 of \cite{Re}). The methods of \cite{VSC}
allows to get very strong results on the asymptotics of iterates
$\mu$ if $\mu$ is symmetric and has a bounded density with compact
support. By combining our methods with the results of \cite{Ba},
\cite{GKR} and \cite{VSC}, we can show that $H$-recurrent groups
have quadratic growth. More precisely we have

\bp\label{hr}
If $G$ is a compactly generated, locally compact second countable $H$-recurrent
group, then there exists closed normal subgroups $H$ and
$K$ of $G$ such that $K$ is compact, $G/H$ is finite and $H/K$ is isomorphic to
a closed subgroup of $\mG _2$.  In particular, any locally compact $H$-recurrent
group has quadratic growth.
\ep

We need the following

\bl\label{lhr}
Assume that $G$ is as in Proposition \ref{hr}.  Then there exists a compactly
supported recurrent $\nu\in M^1(G)$ such that $\nu$ is symmetric with a
bounded density that is bounded from below on an open neighborhood of $e$
generating $G$.
\el

\bo
Using Lemma 2.2 of \cite{Ba} and replacing
$\mu$ by ${1\over 2}(\mu +\check \mu)$, we may assume that $\mu$ is symmetric.
Since $\mu$ is H-recurrent, we can write for some $k>0$ and
$\rho, \sigma \in M^1(G)$, $\mu ^k= r\rho +(1-r) \sigma$
for some $r \in (0, 1]$ and $\rho$ has symmetric density.  Then, if
$\mu ' =  \sum _0 ^\infty {\mu ^n \over 2^{n+1}}$ we have for some $\epsi >0$
$\mu ' \geq {\epsi \over 3}(\rho +\rho \sigma +\sigma \rho)= \epsi\nu '$
where $\nu ' = {1 \over 3}(\rho +\rho \sigma +\sigma \rho)$.
Also using Lemma 2.3 of \cite{Ba} we know that $\mu'$ and $\mu ''
=  \sum _0 ^\infty {(\mu ') ^n \over 2^{n+1}}$ are H-recurrent.  From, above we
get $\mu'' \geq  \sum _0 ^\infty {\epsi ^n (\nu ') ^n \over 2^{n+1}}= \nu ''$
say.  Clearly $\nu ''$ has a symmetric density and since the support of $\nu '$
generates $G$, we get that the support of $\nu ''$ is $G$.  We can in the
above inequality replace $\nu ''$ by $\nu'''$ with a bounded symmetric
density
and the support of $\nu'''$ is $G$.  Then $(\mu '')^2 \geq (\nu ''')^2$ and
$(\nu ''')^2$ has positive continuous density.  Then we can restrict
$(\nu ''')^2$ to a compact neighborhood of identity which generates $G$.
Using Lemma 2.3 of \cite{Ba}, the normalization of $(\nu ''')^2$ is recurrent,
hence H-recurrent.
\eo

In order to make explicit the property of quadratic growth we prove the
following.

\bp\label{qg}
Let $G$ be a compactly generated group of at most quadratic growth.  Then
there exists closed normal subgroups $H$ and $K$ of $G$ such that $K$ is
compact, $G/H$ is finite and $H/K$ is isomorphic to a closed subgroup of
$\mG _2$, the motion group of the plane.
\ep

\bo
By \cite{Lo}, we can assume that $G$ is a real Lie group.  Let $G_0$ be the
component of the identity.  Then $G_0$ is open and $G/G_0$ also has at
most quadratic growth.  Since $G/G_0$ is finitely generated, $G/G_0$
contains a subgroup of finite index isomorphic to $\Z ^i$ for $i\leq 2$.
So, we may assume that $G/G_0 \simeq \Z ^i$ for $i\leq 2$.

Now, $G_0$ has at most quadratic growth implies that $G_0 \simeq K\ltimes
\R^j$ for $j \leq 2$ and $K$ is a compact group.  Then there exists a
characteristic compact subgroup $L$ of $G_0$ contained in $K$ such that
$K/L$ is a subgroup of $SO(2, \R)$.  Replacing $G$ by $G/L$, we may assume
that $K$ is a subgroup of $SO(2, \R)$.

If $j =1$, $G_0 = K\times \R$ and $K$ is characteristic in $G_0$, hence
normal in $G$.  Replacing $G$ by $G/K$ we may assume that $G_0 = \R$.
Then since $G$ is unimodular $G_0$ is central in a subgroup of finite index
in $G$.  Thus, we may assume that $G_0$ is central in $G$, hence
$G$ is nilpotent.  Since $G$ has at most quadratic growth, $G$ is
abelian, hence is a closed subgroup of $\R^2$.

Let $j =2$ and $a \in G$ but $a \not \in G_0$. Then the closed subgroup
$A$ generated by $a$ is isomorphic to $\Z$ and $A\cap G_0 = \{ e \}$
as $G/G_0 \simeq \Z ^i$.  Let $H_a$ be the closed subgroup generated
by $G_0 $ and $a$.  Then $H_a = A\ltimes G_0$ where the action of
$A$ on $G_0$ is given by the conjugacy.  Since $A\simeq \Z$, an easy
computation shows that $H_a$ has growth at least three. This is a
contradiction to $G$ having quadratic growth.  Thus, $G= G_0$ which
is isomorphic to a closed subgroup of $\mG _2$. \eo

\bo $\!\!\!\!\!$ {\bf of Proposition \ref{hr}}\ \
Let $G$ be $H$-recurrent and $L$ be any compactly generated open subgroup of
$G$.  Let $\mu \in M^1(G)$ be $H$-recurrent.  Let $\mu _L \in M^1(L)$ be the
induced measure on $L$.  Then as observed in section 2 $\mu _L$ is also
$H$-recurrent.  Hence we can assume $\mu _L\in M^1(L)$ is as in the
conclusion of Lemma \ref{lhr}.  Then since $L$ is unimodular, we are in the
situation of Theorem VII 1.1 of \cite{VSC}.  It follows that $L$ and hence
$G$ has polynomial growth of degree at most two and the rest follows from
Proposition \ref{qg}.
\eo

\end{section}

\begin{section}{Examples of recurrent or transient behaviors for random walks
on $G$-spaces}

Here $E= G/H$ will be a homogeneous space and $\mu \in M^1(G)$ defines a random
walk on $E$ with trajectories $X_n(\oa )x$ ($x\in E$).  We recall that the
associated Markov operator $P$ satisfies:
$$\sum _0 ^\infty P^k \psi (x) = \int \sum _0 ^\infty \psi (X_n(\oa )x) d\mP (\oa
)$$ for all $\psi \in C_b^+ (E)$.  If $E=G$, there is a natural
$P$-invariant measure, that is left Haar measure.  Here in various cases the
discussion involves a natural $P$-invariant measure.
For a given $x \in E$, we define the following properties that may
or may not be satisfied.

$R_x$: There exists a compact set $K_x \subset E$ such that $\mP $-a.e,
$X_n (\oa )x \in K_x$ infinitely often.

$T_x$: For any compact set $K\subset E$, $\mP$-a.e, there exists $n(\omega ) \in
\N$ with $X_n(\omega ) x \not \in K$ for $n \geq n(\omega )$.

We observe that if property $R$ (defined in 2.3) is valid, then since $E$ is a
countable union of compact sets, property $R_x$ is valid a.e.  We will 
also consider the following property $R^a$ for $(E, P)$:

$R^a$: There exists a compact set $K \subset E$ such that for each $x \in
E$, $\mP $-a.e, $X_n (\oa )x \in K$ infinitely often.

It is easy to see that property $R^a$ implies $R_x$ for all $x\in E$
with $K_x = K$ independent of $x \in E$. Clearly if $R_x$ is valid,
then $\sum _0 ^\infty P^n 1_{K_x} (x) = +\infty $. Also, if $\sum _0
^\infty P^n 1_{K} (x) < +\infty $ for any compact $K$, then $T_x$ is
valid.  If $G_\mu$ is "large" one can expect for "most" $x\in E$,
the trajectories $X_n (\omega ) x$ to have similar asymptotic
behaviors.  More precisely one can expect the existence of a large
set of points $x \in E$ such that $T_x$ (or its complement) is valid; see
\cite{HR} for a discussion of analogous properties when $\mu$ has
density.  In the examples studied below, there exists a natural
$P$-stationary measure $\eta $ on $E$ and we also discuss
recurrence and ergodicity properties with respect to $\eta$: see \cite{Fu1} for 
a discussion on this problem and 
\cite{Pe} for other examples of probabilistic significance.  Property $R$ is 
valid with respect to $\eta$ in examples 1-5, $R_x$ is valid for any $x$ in 
examples 1,2,4, while in example 6, $T_x$ is valid $\eta $-a.e.

(1) {\bf Homogeneous spaces with finite stationary measure:} If $E=
G/H$ is compact, Markov-Kakutani theorem implies the existence of a
$\mu$-stationary probability $\eta$, that is, $\mu *\eta = \eta$. If
$\eta$ is extremal, then ergodicity of $(\Oa \times E, \tilde \theta
, \mP \otimes \eta)$, hence of its natural extension $(\hat \Oa
\times $E$, \hat \theta , \hat \eta )$ is valid.  This situation
takes place also for some noncompact spaces.  For example, if $\mF$
is a local field we can take $G = GL(d, \mF) \ltimes\mF ^d$, the
affine group of $\mF ^d$, and we denote $g \in G$ as $g = (a(g),
b(g))$ with $a(g) \in GL(d, \mF)$ and $b(g) \in \mF ^d$.  Then if
$\mu \in M^1(G)$ and
$$\int (|\log (||a(g)||)|+|\log (||b(g)||)|) d\mu (g) < +\infty $$
$$\lim _{n\to \infty} {1\over n} \int \log ||a(g)|| d\mu ^n(g) <0$$
there exists a unique $\mu$-stationary probability $\eta$
on $\mF ^d= E= G/GL(d, \mF)$ (cf. \cite{CKW}, \cite{Ke2}).  If the support
of $\mu$ has no fixed point on $\mF ^d$, then $\eta$ is not a Dirac
measure. In this case $(E, \eta)$ is a $\mu$-boundary of $(G, \mu)$
(see \cite{Fu}) and remarkable homogeneity properties of $\eta$ at
infinity have been described in \cite{Ke2}.  Using proximality of
the $G_\mu$-action, it is easy to show that property $R^a$ is valid
for $(E, P)$ in this case for any compact subset $K$ having
non-empty interior with $\eta (K)>0$.

Another kind of example takes place when $G$ is a semisimple noncompact
real Lie group and $E$
is a non-compact homogeneous space of finite volume $E= G/\Ga$ where $\Ga$ is a
lattice in $G$.  Then if $\mu \in M^1(G)$ is such that the support of $\mu$ is
compact and generates a Zariski-dense subgroup, then
by \cite{EM} there exists a compact
$K\subset E$ and $C_K >0$ such that for any $x \in E$, $$\liminf _{n \to
\infty}\mu ^n *\delta _x(K) \geq C_K.$$  It follows that if $u \in C_C^+(E)$
satisfies $u\geq 1$ on $K$, then $\sum _0 ^\infty P^ku = +\infty$ on $E$, hence
property $R$ is valid with respect to any $\mu$-stationary measure.
If $G$ is almost simple and for any $g$ in $G$ the group $gG _\mu g^{-1}$ has no 
finite index subgroup contained in $\Gamma$, then it follows 
from \cite{BQ} that Haar measure on $E$ is the unique $\mu$-stationary measure. 
Then, by Breiman's law of large numbers for Markov chains, property $R^a$ is 
valid for any compact set $K$ with positive measure.  Finally, Haar measure is 
$\mu$-stationary but, in general, there exists finitely supported 
$\mu$-stationary measures.

(2) {\bf Affine space and affine groups:} If $\mF$ is a local field with
absolute value $|\cdot |$ and $E$ is
the corresponding affine line, that is,
$G$ is the affine group $'ax+b'$ of $\mF$, then the conditions
$$\int (|\log |a(g)| |+ |\log |b(g)||)^{2+\delta }  d\mu (g) < +\infty ~~(\delta
>0) $$
$$\int \log  |a(g)| d\mu (g) = 0$$
and the support of $\mu$ has no fixed point on $E$ imply the
existence of $u \in C_c^+(\mF)$ such that $\sum _0 ^\infty P^k u =
+\infty$ on $\mF$ (see Babillot and others \cite{BBE}), hence by
Proposition \ref{p1}, the existence of $\lam$ with $P\lam = \lam$.  As
shown in \cite{BBE}, the Radon measure $\lam$ has infinite mass and
is unique up to a coefficient.  Hence the equidistribution property
of Proposition \ref{p1}  is valid in this case.  Here we complete
the example with the following.

\bp\label{ba}
With the above notations and hypothesis, there exists a unique 
$\mu$-stationary Radon measure $\lam$ on $E$, the mass of $\lam$ is infinite, 
hence the system $(\hat \Oa \times E, \hat \theta , \hat \lam )$ is ergodic.  
Furthermore, property $R^a$ is valid.
\ep

\bo
The proof is based on the method of \cite{BBE}.  Let $\tau$ be the first
descending ladder index of the random walk $\log |a(X_n)|\in \R$, that is ,
$\tau (\oa ) = \inf \{ n \in \N \mid |a(X_n(\oa ))| <1 \}$.  Let $\mu _\tau$
be the law of $X_\tau (\oa ) \in G$ and observe that $m_\tau =\int \log |a(g)|
d\mu _\tau (g) <0$.  Also, using the fact that for $t>0$, $$\mP \{ \tau >t
\} \leq ct^{-{1\over 2}}$$ with $c>0$ which follows from fluctuation theory
of random walks on $\R$ (cf. \cite{CKW}),  we get that
$$\mE (|\log (|a(X_\tau )|)|+|\log (|b(X_\tau )|)|)<+\infty .$$
Then we know from the above examples that there exists a
unique $\mu _\tau$-stationary probability $\nu _\tau$ on $E$.  We denote
by $X_{\tau ,n}(\oa )$ the random walk on $G$ defined by $\mu _\tau$, by
$P^\tau$ the corresponding Markov operator on $E$.  Then $\tau$ is a stopping
time and $X_{\tau ,n}(\oa)$ is a subprocess of $X_n(\oa )$.  Since property
$R^a$ is valid for $(E, P^\tau)$, it is also valid for $(E, P)$.

Let $U$ be a relatively compact open subset of $E$ with $\nu _\tau(U)>0$.
Then, from the examples
considered above, we know that, for any $x \in U$, $\mP$-a.e
$X_{\tau , n}(\oa )x \in U$ for some $n=n(\oa ,x) \in \N$.
Hence, also $X_n(\oa ) x\in U$ for some $n \in \N$.
Using the remark following definition of property $R$, we get that
property $R$ is valid for $(E, P, \lam )$.  The uniqueness property of
$P$-invariant Radon measures implies the extremality of $\lam$.  Hence ergodicity
of $(\hat \Oa \times E , \hat \theta , \hat \lam )$ follows from
Corollary \ref{co3}.
\eo

(3) {\bf Fibered random walks:}  We consider below two different
situations with $E= G/\Da '$ and $\Da'$ is a normal subgroup of a
closed subgroup $\Da$ of $G$ such that $G/\Da = \ol E$ is compact
and $\Da /\Da '$ is isomorphic to $Z= \R$ or $\Z$.  Here we give a
general condition of recurrence for these choices which will be
applied below in two different situations.  Since $\Da'$ is normal
in $\Da$, $Z= \Da /\Da '$ acts on the right on $E= G/\Da '$, this
action commutes with the $G$-action and the corresponding factor
space of $E$ is $\ol E$.  The action of $z \in Z$ on $y \in G/\Da '$
will be denote by $y\ct z$.  Also we denote by $\ol P$ the
convolution operator on $\ol E$ defined by $\mu$, and we consider
$\rho \in M^1(\ol E)$ with $\ol P \rho = \rho $ and $\rho$ is
extremal with respect to this property.  Then we consider the space
$\Omega \times \ol E$ and the map $\ol \theta $ with $\ol \theta
(\oa, x) = (\theta \oa , Y_1(\oa ) x)$ if $(\oa, x) \in \Oa \times
\ol E$ and $\theta $ is the shift on $\Oa$.  We endow $\Oa \times
\ol E$ with the measure $\mP \otimes \rho$ and we observe that $\mP
\otimes \rho$ is $\ol \theta $-invariant and ergodic.  We fix a
Borel fundamental domain $D$ of $Z$-action on $E$ which is
relatively compact and $E$ is Borel isomorphic to $D\times Z$.  We
denote by $l$ the measure on $Z$ which is Lebesgue if $Z= \R$ or
counting if $Z = \Z$.  Then we can identify $\rho \otimes l$ with a
Radon measure $\lam$ on $E$ which satisfies $P\lam = \lam$. We define a
Borel function $z(y)$ on $E$ by $y = \ol y \ct z(y)$ where $\ol y
\in D$, $z(y) \in Z$. Also we write $z(g, \ol x) = z(gx)$ if $g \in
G$ and $x \in D$ corresponds to $\ol x \in \ol E$. Then we have the
cocycle relation
$$z(gh, \ol x) = z(g, h\ol x)+z(h, \ol x)$$ for $g , h \in G$ and $\ol x
\in \ol E$.  It follows that for $g _i \in G$ ($1\leq i \leq n$) and $\ol
x \in \ol E$,
$$z(g_n \cdots g_1, \ol x)=\sum _1 ^n z(g_k , g_{k-1}\cdots g_1\ol x).$$
If we write $f(\oa, \ol x) = z(Y_1(\oa) ,\ol x )$, $S_n (\oa, \ol x
) = z(Y_n \cdots Y_1, \ol x)$ with $\omega \in \Oa $ and $\ol x \in
\ol E$, we get $S_n(\oa , \ol x) = \sum _0 ^{n-1} f( \theta ^{k}
(\omega, x))$. With the above identification, we write $x = (\ol
x, t) \in \ol E \times Z$ and $X_n(\oa ) x= (X_n(\oa )\ol x ,
t+S_n(\oa ,\ol x))$ where $X_n(\oa )\ol x$ is the random walk on
$\ol E$ defined by $\mu$. We will also consider the map $\tilde
\theta $ on $\Oa \times E$ defined by $\tilde \theta (\oa, y) =
(\theta \oa , Y_1 (\oa )y)$ and $\Oa \times E$ will be endowed with
the $\tilde \theta$-invariant measure $\mP \otimes \lam$. This is a
skew-product as in \cite{Gu3} and \cite{Sc} with the base $\Oa
\times \ol E$, the fiber $Z$ and the function $f(\oa , x)$.  We
refer to \cite{Gu3} for a discussion of ergodic properties of such
systems in case $\rho$ is $G$-invariant.
Since $D$ is relatively compact for $g \in G$, $c(g) =
\sup _{x\in D} |z(gx)| < +\infty$.  Also $c(gh) \leq c(g)+c(h)$.  In
particular, if $G$ is a closed subgroup of a linear group, it
follows that for some $c>0$, $c(g) \leq c(\log ||g|| +1)$. Also it
is easy to show that if $(D, c)$ is replaced by $(D', c')$ for
another fundamental domain $D'$, then there exist constants $d,
d'\geq 0$ such that $c'(g) \leq dc(g)+d'$.  Thus, it follows for
$\mu \in M^1(G)$ that the condition $\int c(g) d\mu (g) < +\infty$
is independent of $(D, c)$. It is also evident that if $\int c(g)
d\mu (g) < +\infty$, then we get the $\mP \otimes
\rho$-integrability of $f(\oa, \ol x)$, hence we can apply ergodic
theorems to the Birkhoff sum $S_n(\oa , \ol x)$.  Also in this case we
can describe the measure $\hat \lam$ on $\hat \Oa \times E$ in terms of
$\lim _{n\to \infty} Y_{-n}(\hat \oa ) \cdots Y_0(\hat \oa )\rho= \rho _
{\hat \oa }$ which is the limit of the bounded martingale
$Y_{-n}(\hat \oa ) \cdots Y_0(\hat \oa )\rho$.

\bp\label{r0} With the above notations, assume that $\mu \in M^1(G)$
is such that $\int c(g) d\mu (g) < +\infty$. Then, if $\int z(gx)
d\mu (g) d\rho (\ol x) = 0$, the random walk on $G/\Da '$ defined by
$\mu$ satisfies property $R$ with respect to $\lam = \rho \otimes l$.  
We have $\hat \lam = \int \delta _{\hat \oa }\otimes l _{\hat \oa} 
d\hat \mP (\hat \oa )$, where $\rho _{\hat \oa } = \lim _{n\to \infty} 
Y_{-n}(\hat \oa )\cdots Y_0(\hat \oa ) \rho $ and 
$l_{\hat \oa} = \rho _{\hat \oa }\otimes l$.  If $\rho$ is 
$G_\mu$-invariant
and $\mu$ is symmetric, then the condition $$\int z(gx) d\mu (g)
d\rho (x) =0$$ is satisfied. \ep

\bo Let $I\subset Z$ be an open symmetric interval. We write $x =
(\ol x, t)$ with $\ol x \in D$, $t \in I$, $X_n(\oa )x = (X_n(\oa
)\ol x, t +S_n(\oa ,\ol x))$. Since $\int f(\oa, \ol x) d\mP (\oa )
d\rho (\ol x) =0$, the Birkhoff sum $S_n(\oa, x)$ belongs to $I$
infinitely often $\mP \otimes \rho$-a.e (see \cite{Sc}).  It follows
that for any $t\in I$ and $\mP \otimes \rho$-a.e $X_n(\oa )x \in
D\times I$ infinitely often.  Since the actions of $G$ and $\Da /\Da
'$ commute, the same property is true for $s \in Z$ with $x$ 
replaced by $x\cdot s$ and $I$ by $I+s$.  Then using Poincar\'e
recurrence theorem, for $D\times 1_I$ and $\rho \otimes (1_Il)$,
property $R$ with respect to $\lam$ follows.  We observe that the
structure of skew product implies $g(\rho \otimes l) = (g\rho )
\otimes l$ if $g \in G$, hence for any $\psi \in C_c(E)$,
$Y_{-n}(\hat \oa ) \cdots Y_0(\hat \oa)\lam (\psi) = (Y_{-n}(\hat
\oa )\cdots Y_0(\hat \oa )\rho) \otimes l(\psi)$ is a bounded
martingale which converges to $\rho _{\hat \oa} \otimes l (\psi)$.
Hence from remark \ref{nrk} we get that $\hat \lam = \int \delta 
_{\hat \oa } \otimes l_{\hat \oa } d\hat\mP (\hat \oa)$.

The cocycle property
of $z(g, \ol x)$ gives $z(g, \ol x) = -z(g^{-1}, g\ol x)$ on $G\times \ol E$.
Then $\int z(g, \ol x ) d\mu (g) d\rho (\ol x) = - \int z(g^{-1}, g\ol x)
 d\mu (g) d\rho (\ol x)$.  The $G_\mu$-invariance of $\rho$ and the symmetry of
$\mu$ gives $\int z(g, \ol x ) d\mu (g) d\rho (\ol x) = - \int z(g, \ol x)
 d\mu (g) d\rho (\ol x)$.  Thus proving the last assertion.
\eo

(4) {\bf Pointed vector spaces:}  Let $G= GL(d, \R)$ with $d \geq 2$.  We
assume
that $\mu \in M^1(G)$ satisfies the following conditions (H-1)-(H-3).

\be
\item[(H-1)] no finite union of proper subspaces of $\R ^d$ is invariant
under the support of $\mu$, that is, $G_\mu$ is strongly irreducible on
$\R ^d$.

\item [(H-2)] $G_\mu$ contains an element $g$ that has unique dominant
eigenvalue $\lam _g$ with $|\lam _g| = \lim _{n \to \infty}
||g^n ||^{1\over n}$.

\item [(H-3)] the support of $\mu$ is compact.
\ee

Here the space $E$ will be the factor space $V= (\R ^d )^*/ \{\pm Id \}$
considered
as a $G$-homogeneous space where $(\R ^d )^* = \{ v\in \R ^d \mid v\not = 0 \}$.
We denote by $e$ the
point of $V$ corresponding to the first basis vector, by $\Da '$ the
stabilizer of $e$ and $\Da$ the stabilizer of the line $\R _+^* \subset
V$ containing $e$.  Then $\Da'$ is normal in $\Da$ with $\Da / \Da' \simeq \R _+^*$ and
$E=V= G/\Da '$, $\ol E = G/ \Da \simeq \mP ^{d-1}$, the projective space.
Also $v\in V$ can be written as $v = (x, r)$ with $r \in \R _+^*$ and $x \in
\mP ^{d-1}$ so that $V= \mP ^{d-1} \times \R _+^*$ and $\mP^{d-1}$ can be
considered as a fundamental domain of the action of $\R _+^* \simeq \Da/\Da
'$ on $V$ by dilations.  We identify $\mP ^{d-1} \subset V$ with the unit
sphere modulo symmetry.  In particular, if $g \in G$ and $x \in V$, the
quantity $||gx||$ is well defined and gives a cocycle on $G\times \mP
^{d-1}$.  Also, if $x, x' \in \mP ^{d-1}$, then
$\ol \delta (x, x') = |\sin (x,x')|$
defines a distance on $\mP ^{d-1}$.  We denote by $H_\epsilon$
the space of $\epsilon$-Holder functions on $\mP ^{d-1}$ endowed with the
norm
$$||\phi ||_\epsilon = \sup _{x\in \mP^{d-1}} |\phi (x)| + \sup _{x, x'\in
\mP^{d-1}} {|\phi (x) -\phi (x')| \over {\ol \delta } ^\epsilon (x, x')}.$$
Then under condition (H) for $\epsilon $ sufficiently small, the operator
$\ol P$ on $\mP ^{d-1}$ associated with $\mu$ has nice spectral
properties on $H_\epsi$ (cf. \cite{GuS}).  In particular, $\ol P$ has a unique
stationary measure $\rho = \nu$ on $\mP ^{d-1}$, and $\ol P$ has a spectral gap
on $H_\epsilon $, i.e $\ol P$ is the sum of a one-dimensional projection $\pi$
and another operator $Q$ with spectral radius less than one which satisfies
$Q\pi = \pi Q =0$.  We need to consider also the family of operators $P_t $ ($t
\in \R$) defined by $ P_t \phi (x) = \int ||gx||^{it} \phi (gx) d\mu (g)$;
we refer to \cite{GuS} for a recent exposition of the spectral
properties of the operators $ P_t$.  In particular, if $t \not = 0$, $P_t $ has a
spectral radius less than $1$, and for $t$ small one has
$ P_t = k(t ) \pi _t + R_t$ where $k(t ) \in \C$, $\pi _t$ is a projection of
rank one, $R_t$ has  a spectral radius less
than $|k(t )|$ and commutes with $\pi _t$.  All quantities depend
analytically of $t$, and if $\int \log ||gx|| d\mu (g)
d\nu (x) = 0$, then $k '(0) = 0$ and $\sigma ^2 = -k''(0) >0$ (see \cite{BL}).
We denote by $l$ the Lebesgue measure ${dr\over r}$ on $\R ^*_+$, we
observe that the Radon measure $\lam =\nu \otimes l$ on $E= V$ is
$P$-invariant and we have the following.

\bp\label{r1}
Assume $\mu \in M^1(G)$ satisfies conditions (H-1)-(H-3) and $\int \log ||gx||
d\mu (g) d\nu (x) = 0$.  Then there exists a $\sigma >0$ such that for any
$\psi \in C_c(V)$ and any $v \in V = (\R ^d )^*/\{\pm Id\}$,
$$\lim _{n \to \infty} \sigma \sqrt { 2\pi n} P^n \psi (v) = (\nu \otimes
l )(\psi).$$ In particular, for any $u \in C_c^+(V)$ with $(\nu
\otimes l )(u) >0$, $\sum _0 ^\infty P^k u = +\infty$ on $V$.
Furthermore property $R^a$ is valid , $\nu\otimes l$ is
$P$-invariant extremal, hence $(\hat \Oa \times V , \hat \theta ,
\hat \lam )$ is ergodic. \ep

\bo
A general formula as in the proposition was proved in \cite{LeP}, under the
assumption that $P$ has spectral gap and the spectral radius of $P_t$ ($t\not
=0$) is less than one, which is valid in our situation (cf. \cite{GuS}).  Hence,
here we only sketch the proof.  For fixed $v \in V$, we consider the sequence
of Radon measures $\mu _n$ defined by $\mu _n (\psi ) = \sigma \sqrt { 2 \pi n}
P^n \psi (v)$.  It suffices to prove the convergence of $\mu _n$ to $\nu
\otimes l (\psi)$ for functions of the form $\psi (x, r) = \phi (x)
f(r)$ where $\phi \in H_\epsi$ and $f \in L^1(\R _+^*)$ is such that its
Fourier transform $\hat f$ has compact support.  Then we can write,
using Fourier inversion formula $f(r) = {1\over 2\pi} \int \hat f (t )
r^{-it } dt$ where $\hat f (t ) = \int _{R_+^*} r^{it } f(r)
{dr \over r}$.  Then $P^n \psi (x, r) = {1\over 2\pi} \int P_t ^n \phi
(x) r^{-it } \hat f (t ) dt$.  Replacing $t$ by $t \over
\sqrt n $ we get that $$\mu _n (\psi ) = {\sigma \over \sqrt {2\pi }}
\int P_{t \over \sqrt n} ^n \phi (x) r^{-i{t \over \sqrt n}} \hat f
({t \over \sqrt n}) dt .$$  Using the spectral properties of
$P_t$ we can replace $P_{t \over \sqrt n} \phi (x)$ by $k ^n ({t
\over \sqrt n}) \pi _{t \over \sqrt n} \phi$ which converges to
$e^{-\sigma ^2 t ^2 \over 2} \nu (\phi )$.  Then (see \cite{LeP}) we
get that
$$\lim \mu _n (\psi ) = \nu (\phi ) \hat f(0) {\sigma \over
\sqrt{2\pi }} \int
e^{-\sigma ^2 t ^2 \over 2} dt = (\nu \otimes l)(\psi).$$

In order to show property $R^a$, we use a result of \cite{Bo}.  Under
conditions (H), for any $v \in V$ we have $\mP$-a.e
$$\limsup  _{n\to \infty}||X_n(\oa ) v|| = \infty, ~~~
\liminf _{n\to \infty} ||X_n(\oa )v|| =0.$$  Since the support of $\mu$ is
compact, there exists $c>0$ such that for any $(\oa ,v)\in \Oa \times V$, $n
\in \N$,
$${1\over c} \leq {||X_{n+1}(\oa )v||\over ||X_n(\oa )v||} < c.$$  Then, for
any $v$, the relatively compact open set $$U_c = \{ u\in V\mid {1\over c}<
||u|| <c \}$$ is visited infinitely often $\mP$-a.e.  Hence property $R^a$
is valid.  The validity of property $R$ follows from Proposition
\ref{r0}.

The extremality of $\lam = \nu \otimes l$ is proved in Proposition 3.5 of
\cite{GuS}, hence the ergodicity of $(\hat \Oa \times V, \hat \theta , \hat \lam)$
follows from Corollary \ref{co3}.
\eo

\br In general, the support of $\nu$ in $\mP ^{d-1}$ is of Cantor
type, hence the support of $\nu \otimes l $ is a proper subset of
$V$.  If $u \in C_c^+(V) $ with $(\nu \otimes l )(u) =0$, we can
show that $\mP$-a.e, for any $v$, $\sum _0 ^\infty u(X_n(\oa )v)
<+\infty$, in particular, if $K$ is a compact set that does not
intersect the support of $\nu \otimes l$, then $\mP$-a.e the random
walk $X_n(\oa )v$ escapes from $K$, if $v$ is not in the support of
$\nu \otimes l$. \er

(5) {\bf Covering spaces:} Let $G$ be a real simple linear group of
rank one, $\Da$ be a cocompact lattice and $\Da '$ be a normal
subgroup in $\Da$ such that $\Da /\Da ' \simeq \Z$.  It may be noted
that if $G$ has rank more than one, then $G$ has Kazdhan property T
which implies that such a $\Da '$ does not exist. Let $E= G/\Da '$,
$\ol E = G /\Da $ and $D$ be a fundamental domain of the $\Z$-action
on $E$ as in (3).  Let $m$ (resp. $\ol m$) denote the Haar measure
on $E$ (resp. $\ol E$), hence with the notations of (3) $\rho = \ol
m$ and $m =\rho \otimes l$.  Such a situation arises if $G= SL(2,
\R)$ and $\Da$ is the fundamental group of a compact Riemann surface
$S$ and in this case $G/\Da$ (resp. $G/\Da '$) can be identified
with the tangent unit bundle of $S$ (resp. of an abelian cover of
$S$).

\bp\label{r2}
Assume $\mu \in M^1(G)$ is symmetric with $\int \log ||g|| d\mu (g) < +\infty$
and $G_\mu$ is non-amenable.  Then $(\hat \Oa \times E, \hat \theta , \hat
\mP \otimes m)$ is ergodic.
\ep

The proof depends on the next two lemmas.

\bl\label{rl1}
Suppose $\mu \in M^1(G)$ is such that $G_\mu$ is non-amenable.  Let $\ap$ be a
non-trivial character of $\Da$ and $T _\ap$ be the unitary representation of
$G$ induced by $\ap$.  Then the spectral radius of $T _\ap (\mu )$ is less
than one.
\el

\bo
In view of Theorem C of \cite{Sh}, it suffices to verify that $T_\ap$
does not weakly contain the trivial one-dimensional representation of $G$.
It is clear that the one-dimensional representation of $\Da$ defined
by $\ap$ does not weakly contain the trivial one-dimensional
representation of $\Da$.  Since $\Da$ is a cocompact lattice in $G$, it follows
from 1.10 and 1.11, Chapter III of \cite{Ma} that the induced representation
$T_\ap$ does not weakly contain the trivial one-dimensional representation of
$G$.
\eo

If $u \in L^2 (\ol E, \ol m)$, $f \in l^2(\Z)$, we can identify
$u \otimes f$ with an element of $L^2(E, m)$, again denoted by
$u\otimes f$.  If $f \in l^1 (\Z )$ with $\sum f(k) = 0$, we write
$f \in l^1_0 (\Z )$.  In Proposition 3.6 of \cite{GuS}, if $X$ is a compact metric
space, a Markov operator
$Q$ acting on $Y = X\times \R ^d$ which commutes
with the $\R ^d$ translations was considered and it was proved that $\lim _{n
\to \infty } ||Q^n (u\otimes f) || _1 = 0$ for Holder functions $u \in H_\epsi
(X)$ and $f \in L^1(\R ^d)$ with $\int f(x) dx =0$.  Essential points of the
proof are polynomial growth of $\R ^d$ and a spectral gap
property for the
$Q$-action on functions of the form $u \otimes \ap$ where $u \in H_\epsi (X)$
and $\ap$ is a given character of $\R ^d$ (see above).

Here we observe that the adjoint $P^*$ of $P$ in $L^2(E, m)$ is associated with
$\check \mu$, the symmetric of $\mu$ which has the same properties as $\mu$.
Also, in view of Lemma \ref{rl1}, the condition of spectral gap of $P^*$
is valid for the $P^*$-action on functions of the form $u \otimes \ap$ where
$u \in L^2(\ol E, \ol m )$ and $\ap$ is a fixed character of $\Z$.  The
polynomial growth condition is also satisfied here.  Hence the proof in
\cite{GuS} gives the following.

\bl\label{rl2}
Assume $\mu$ satisfies the condition in Lemma \ref{rl1}.  Then for any $u \in
L^2(\ol E, \ol m )$ and $f \in l^1_0( \Z)$, we have $\lim _{n \ra \infty}
||{P^*} ^n (u\otimes f) || _1 =0$.
\el

\bo $\!\!\!\!\!$ {\bf of Proposition \ref{r2}}\ \
We begin by showing that if $h \in L^\infty (E, m)$ satisfies $Ph = h$, then $h$
is constant $m$-a.e.  If $u \in L^2 (\ol E, \ol m )$ and $f \in l^1 _0 (\Z)$,
then $<(P^*)^n (u \otimes f), h> = <u\otimes f, h>$.  By Lemma \ref{rl2},
$<u \otimes f, h> = 0$.  Since $f$ is arbitrary in $l^1 _0(\Z)$ this
implies the $\Z$-invariance of $h$.  Hence $h$ defines an element
$\ol h \in L^\infty ( \ol E, \ol m)$ such that $\ol P (\ol h )= \ol h$.
This equation can be written as $$ \ol P (\ol h )(\ol x)=
\int \ol h (g \ol x) d\mu (g) = \ol h (\ol x).$$  Since $\ol m$ is
$G$-invariant and $\ol h \in L^2 (\ol E, \ol m)$, strict convexity in $L^2(\ol E,
\ol m)$ implies that $\ol h(\ol x) = \ol h (g \ol x)$ $\mu \otimes \ol m$-a.e.
Since $G_\mu$ is non-amenable, it is unbounded, hence the Moore ergodicity
theorem implies the ergodicity of $G_\mu$ on $\ol E = G/\Da $
with respect to $\ol m$ (cf. Theorem 4 of \cite{Mo}), hence $\ol h$ is constant 
$\ol m$-a.e.  This proves that $h$ is constant $m$-a.e.

The above argument shows also the extremality of $\rho = \ol m$ as a $\ol
P$-stationary measure.  The condition $\int \log ||g|| d\mu (g) <+\infty$
implies that the condition $\int c(g) d\mu (g) <+\infty$ of Proposition
\ref{r0} is satisfied.  Then Proposition \ref{r0} gives that property $R$ is
valid with respect to $m = \ol m \otimes l$.  Then as above, we can use
Proposition \ref{ec} to obtain the required ergodicity.
\eo

\br
In the situation of Proposition \ref{r2}, a counter example to property
$R^a$ is the following.  Assume $G_\mu = \Da'$.  Since $\Da$ is
non-amenable, the same is true of $\Da '$ and all the conditions on $\mu$
can be seen to be satisfied.  Since $\Da '$ is normal in $\Da$, we observe
that if $e \in E$ corresponds to $\Da '$ and $z \in \Da$, then $e\cdot z$
is $\Da '$-invariant.  It follows that $P(e\cdot z, \cdot ) =
\delta _{e\cdot z}$ and if $u(e\cdot z)=0$, that is $x=e\cdot z$ is not in the
support of $u$, then $X_n(\oa )x$ is not in the support of $u$ for any $n$
and a.e $\omega $.  We conjecture that property $R_x$ is valid for any $x$ and
if $\mu$ is adapted, then property $R^a$ is valid.
\er

(6) {\bf Transient behavior on pseudo-Riemannian symmetric spaces:}  Let $G$
be a semisimple real algebraic group with no compact factors and $H$ be a
closed subgroup of $G$ that is not Zariski dense in $G$.  Assume that 
$E=G/H$ has a $G$-invariant measure $m$ and $\mu \in M^1 (G)$ is 
such that $G_\mu$ is non-amenable.
Then $P$ acts on $L^2 (E, m)$ as a contraction.  As an extension of the Borel
density theorem it is proved in \cite{Gu4} that, since $H$ is not
Zariski-dense, $G/H$ do not carry an invariant mean, that is, $G/H$ is not 
amenable in the sense of Eymard (see also \cite{St}).  
In particular, the spectral radius of $P$ on $L^2
(E, m)$ is strictly less than one.  Hence for $\phi \in C_c(E)$, one has $\sum _0
^\infty P^k \phi \in L^2 (E, m)$. In particular $\sum _0^\infty P^k \phi $ is
finite $m$-a.e.  This implies property $T_x$ is valid $m$-a.e.

We consider the case where the homogeneous space $G/H$ is a symmetric
pseudo-Riemann space, that is $H$
is the set of fixed points of an involution.  Then $H$ is reductive, hence $H$
is
unimodular.  Thus, $G/H$ has a $G$-invariant measure $m$.  Also, $H$ is
algebraic, hence not Zariski dense.  Then the above discussion gives the
following.

\bp
Assume $E=G/H$ is a pseudo-Riemannian symmetric space, $\mu \in M^1(G)$ is
such that $G_\mu$ is
non-amenable.  Then for $\phi \in C_c(E)$, $\sum _0 ^\infty P^k \phi $ is finite
$m$-a.e.  In particular, $T_x$ is valid $m$-a.e.
\ep

Among the pseudo-Riemannian symmetric spaces we have the spaces $E=
SL(n , \R) /SO(p, q)$ $(p+q < n)$ and $E= SO(k, l)/SO(p, q)$ $(p+q <
k+l)$.  The space $SL(2, \C)/SL(2, \R)$ locally isomorphic to
$SO(3,1)/SO(2,1)$ was considered in \cite{Gu4}.

(7) {\bf On homogeneous spaces with property $R$ :} Considering the
above examples, we formulate the following questions for a $G$-homogeneous
space with $G$ a connected Lie group.

a) Characterize the systems $(E, P, \lam )$ with property $R$.

More precisely, if $R$ is valid is it true that $\lam$ has at most quadratic 
growth in the following
sense:  if $V\subset G$ is a compact neighborhood of $e$, and $x$ in the support
of $\lam$ does there exist $c>0$ such that, for any $n \in \N$,
$\lam (V^nx) \leq cn^2$.  This property is valid in examples 2,4,5 with linear
growth of $\lam (V^nx)$.  It may be noted that when $P$ is given by spread-out 
probability on $G$, question (a) has definitive answer if $G$ is a compact 
extension of simply connected nilpotent real Lie group \cite{GS} or if 
$G$ is a $p$-adic algebraic group \cite{RS}.   

In view of Corollary \ref{co3}, if $\lam$ is extremal, 
property $R$ is equivalent to ergodicity of $(\hat \Oa\times E, \hat \theta , 
\hat \lam )$.

b) Characterize systems $(E, P, \lam)$ with $\lam$ of finite mass.  

Examples of this situation are given in (1).  Except for trivial situations, is 
the general case a combination of these kind of examples.

\end{section}

\begin{section}{Singularity of stationary measures on the projective
line}

As is well known the group $G= SL(2, \R)$ acts by fractional linear
transformations on the Poincar\'e upper half plane $\mH = \{ z \in \C \mid
y = {\rm Im} z >0 \}$, preserving the Poincar\'e metric ${|dz|^2\over y^2}$.
If $a = \pmatrix {1&2 \cr 0 &1 }$ and $b = \pmatrix {1& 0 \cr 2 &1}$, the group
$\Ga = < a,b>$ is a free subgroup of index $6$ in $SL(2, \Z) \subset G$.
Furthermore the $\Ga$-action on $\mH$ is free and totally discontinuous.  We
denote by $\Ga '$ the commutator subgroup of $\Ga$, by $\ol \ga$ the projection
of $\ga \in \Ga$ in $\ol \Ga = \Ga /\Ga ' = < \ol a , \ol b> = \Z ^2$, and
$|\ol \ga|$
the word length of $\ol \ga $ in $\ol a , \ol b$.  We observe that as a Riemann
surface, $\Ga \bs \mH$ can be identified with the complement of $\{ 0, 1\}$ in
$\C$; also $\Ga '\bs \mH$ is
an abelian cover for $\Ga \bs \mH$ with covering group $\Z^2 = <\ol a, \ol b>$.
Then $\Ga \bs G$ (resp. $\Ga '\bs G$) can be identified with the unit tangent
bundle of $\Ga \bs \mH$ (resp. $\Ga ' \bs \mH$).

Let $MAN$ be the triangular subgroup of matrices of the form $\pmatrix {a&b\cr
0&{1\over a}}$ ($a \in \R ^*$ and $b\in \R$).  Here
$M = \{ \pm Id \}$, $A = \{ \pmatrix {a & 0 \cr 0 & {1\over a}} \mid a >0 \}$ and
$N = \{ \pmatrix {1 & b \cr 0 & 1} \mid b \in \R \}$.  Then $G/MAN \simeq \mP ^1$
can be identified with the boundary of $\mH$ in $\C$.  If $\mu \in M^1 (\Ga )$ is
adapted, there exists a unique $\mu$-stationary measure $\nu$ on
$\mP ^1 = G/ MAN$: $\mu *\nu = \sum _{\ga \in \Ga } \mu (\ga ) \ga \nu = \nu$.
It is known (see for instance \cite{Fu} Proposition 4.1, p.207) that Lebesgue
measure $m$ on $\mP ^1$ is
such a measure for a certain $\mu$ as above.  From a geometrical point of view,
such measures appear in the study of foliated Brownian motion on $\Ga \bs G$
along $AN$-orbits and corresponding harmonic measures.  Here, as
an application of recurrence-transience properties of random walk on $\Ga \bs G$
we give a proof of the following (mentioned in \cite{GL}; see \cite{BM} and
\cite{DKN} for recent different proofs).

\bp\label{sp}
Assume that $\mu \in M^1 (\Ga )$ is symmetric adapted and satisfies
$\sum _{\ga \in \Ga} \mu (\ga ) |\ol \ga| ^2 < +\infty$.  Then the unique
$\mu$-stationary measure on $\mP ^1$ is singular with respect to Lebesgue
measure.
\ep

The proof depends on the next two results which are of independent interest.

\bl\label{sl1}
The action of $\Ga '$ on $(\mP ^1 \times \mP ^1, m \otimes m)$ is not ergodic.
\el

\bo
The Lemma is a direct consequence of the results of \cite{Su} and \cite{LM}.
From \cite{LM} we know that Brownian motion on $\Ga '\bs \mH$ is
transient.  Then it follows from Theorem 4 of \cite{Su} that the geodesic flow
on the manifold $\Ga ' \bs \mH$ is not ergodic.
In other words the action of $A$ on $\Ga '\bs G$ endowed with
the Haar measure is not ergodic. By duality this means that the action of $\Ga '$
on $G/A$ endowed with Haar measure is not ergodic.  Now the result follows
since $G/AM$ and $(\mP ^1 \times \mP ^1, m \otimes m)$ are isomorphic as measured
$G$-spaces.
\eo

\bl\label{sl2}
The action of $\Ga '$ on $(\mP ^1 \times \mP ^1, \nu \otimes \nu )$ is ergodic.
\el

\bo
As in 2.3, we denote $\hat \Oa = \Oa ^- \times \Oa$, $\hat \mP = \mP^- \otimes
\mP$,
we write $\hat \Oa = (\oa ^- , \oa )$ if $\hat \oa \in \hat \Oa$.  Then $\theta
\hat \oa = (\oa ^-Y_1(\oa ), \theta \oa )$.  We consider the
transformation
$\hat \theta$ of $\hat \Oa \times \Ga$ defined by $\hat \theta (\hat \oa
,\ga) =
(\theta \hat \oa , \ga Y_1(\oa ))$.  We denote also by $\hat \theta$ the
transformation of $\hat \Oa \times \ol \Ga$ defined by $\hat \theta (\hat
\oa , \ol \ga) = (\theta \hat \oa , \ol \ga \ol {Y_1(\oa )})$.   Since the
random walk on $\Z ^2$ associated with $\mu$ is recurrent and adapted the
skew product $(\hat \Oa \times \ol \Ga , \hat \theta, \hat \mP \otimes \lam _{\ol
\Ga})$ is ergodic (see Corollary \ref{co2}).
If we denote by $z(\oa )$ the random variable on $\mP ^1$ uniquely defined
$\mP$-a.e by $Y_1(\oa) z(\theta \oa ) = z(\oa )$ we have, $\mP$-a.e:
$\delta _{z(\oa )}= \lim _{n} Y_1 \cdots Y _n \nu$
and $\nu$ is the law of $z(\oa )$.  Since $\mu$ is symmetric, $\nu$ is also the
law under $\mP ^-$ of $z'(\oa ^-)$ where $\delta _{z'(\oa ^-)} =
\lim _n Y_0 ^{-1} (\oa )\cdots Y_n^{-1}(\oa )\nu$
and $\mu \otimes \mP ^-$-a.e: $z'(\oa ^-Y_1(\oa ) )= Y_1 ^{-1}(\oa ) z'(\oa ^-)$.
Let $\phi$ be a Borel function on $\mP ^1 \times \mP^1$ and denote
$f( \hat \oa , \ga) = \phi (\ga z(\oa ), \ga z'(\oa ^-))$.  Then the functional
equations above implies
that $f(\theta  \hat \oa , \ga Y_1 (\oa ) ) = f(\hat \oa , \ga)$.  Also,
if $\phi$ is $\Ga '$-invariant $f(\hat \oa, \eta \ga ) = f(\oa ,\ga )$ for all
$\eta \in \Ga '$.  Then we can write $f(\hat \oa , \ga ) = \ol f (\hat \oa ,
\ol \ga )$ where $\ol f$ is a $\hat \theta$-invariant function on
$\hat \Oa \times \ol \Ga $.  From above we know that $\ol f$ is constant
$\mP ^-\otimes \mP$-a.e.  Since the law of
$(z'(\oa ^-), z(\oa ))$ under $\mP ^-\otimes \mP $ is $\nu \otimes \nu$, this
means that
$\phi$ is constant $\nu \otimes \nu$-a.e.  In other words $(\mP ^1 \times \mP^1,
\nu \otimes \nu )$ is $\Ga '$-ergodic.
\eo

\bo $\!\!\!\!\!$ {\bf of Proposition \ref{sp}}\ \
Since $\Ga \bs G$ has finite volume, the geodesic flow on $\Ga \bs G$ is ergodic
with respect to the Haar measure and the action of $A$ on $\Ga \bs G$ is ergodic.
By duality this gives the ergodicity of $\Ga$ on $(\mP ^1, m) = (G/MAN, m)$ .
In particular, any $\Ga$-quasiinvariant and absolutely continuous measure on
$\mP ^1$ is equivalent to $m$.  If $\nu$ is not singular with respect to $m$, we
can write
$\nu = \nu _a + \nu _s$ where $\nu _a \not = 0$ is absolutely continuous and $\nu
_s$ is singular (with respect to $m$).  The equation $\nu = \mu *\nu$ implies $\mu
* \nu _a = \nu _a$ and $\mu *\nu _s = \nu _s$, hence by uniqueness of the
$\mu$-stationary measure, we have $\nu _a = \nu$ and $\nu _s = 0$.  From above,
the $\Ga$-quasiinvariance of $\nu = \nu _a$ implies that $\nu$ is equivalent to
$m$.  Then the two lemmas give the required contradiction.
\eo

\end{section}

\begin{section}{Appendix: Groups whose closed subgroups are unimodular}

Groups whose closed subgroups are unimodular plays a crucial role in
our proof of quadratic growth conjecture.  This motivates us to prove the following
characterization of such groups among algebraic groups and almost connected
locally compact groups in terms of polynomial growth.

\bp\label{ag}
Let $G$ be either the group of $\mF$-rational points of an algebraic group
defined over a non-archimedean local field $\mF$ of characteristic zero or an
almost connected locally compact group.  Then closed subgroups of $G$ are
unimodular if and only if $G$ has polynomial growth.  Furthermore, in the first
case, any compactly generated closed subgroup $H$ of $G$ has a basis of compact
open normal subgroups $(K_n)$ such that $H/K_n$ has a finitely generated abelian
subgroup of finite index.
\ep

\bo
If $G$ is a Lie group over a local field, then Ad$~(g)$ denotes the
adjoint automorphism of the Lie algebra $\cal G$ of $G$ defined by $g$.

We first consider the case of group of $\mF$-rational points of an
algebraic group defined over a non-archimedean local field $\mF$ of
characteristic zero.  Since the Zariski-connected component is a subgroup of
finite index, we may assume that $G$ is Zariski-connected.
Assume that closed subgroups of $G$ are unimodular.  Let $H$ be a open compactly
generated subgroup of $G$.  Then by Corollary \ref{c3}, $H$ contains a
basis $(K_n)$ of compact open normal subgroups.  On the other hand, since Ad$(G)$ 
is also an algebraic group, Ad$(G)$ is closed in $GL({\cal G})$ and is isomorphic 
to $G/Z$ where $Z$ is the center of $G$; as characteristic of $\mF$ is zero, 
center of $G$ is the kernel of the adjoint homomorphism of $G$ (ref. 0.15 and 
0.24 of \cite{Ma}).  Since
$H$ is open in $G$, Ad$(H)$ is open in Ad$(G)$ and hence closed in
$GL({\cal G})$.  Since the orbits of Ad$(H)$ in $\cal G$ are relatively compact,
Ad$(H)$ is also relatively compact, hence Ad$(H)$ is compact.  Since Ad$(G)
\simeq G/Z$, $Z\cap H$ is co-compact in $H$.  Thus, the center of $H$ is
co-compact in $H$.  The same property is valid in the finitely generated subgroup
$H/K_n$.  In particular, $H$ as well as $G$ has polynomial growth.  Now if $H$ is
any compactly generated closed subgroup, let $C$ be a compact generating set
of $H$.  Let $V$ be a compact neighborhood of $C$ and $N$ be the subgroup
generated by $V$.  Then $N$ has a basis $(K_n)$ of compact open normal
subgroup such that $N/K_n$ has a finitely generated abelian subgroup of finite
index.  Let $L_n = K_n \cap H$.  Since $H$ is a closed subgroup of $N$, $L_n$ is
a compact open normal subgroup of $H$ and $H/L_n $ is continuously embedded in
$N/K_n$.  Thus, $H/L_n$ also has a finitely generated abelian subgroup of finite
index.

We now consider a connected real Lie group.  Then for any $g$ in
Ad~$(G)$ (with the quotient topology from $G$), $C_g$ is a simply
connected nilpotent Lie subgroup of $GL({\cal G})$ and $C_g \subset
\{ v \in GL({\cal G}) \mid \lim _{n\to \infty} g^nvg^{-n} = {\rm Id} ~~{\rm in}~~
GL({\cal G}) \} = \tilde C_g$, say.  But $\tilde C_g$ is a unipotent
algebraic group.  Since $C_g$ is a connected Lie subgroup of $\tilde
C_g$ which is unipotent, we get that $C_g$ is closed in Ad~$(G)$. By
Lemma \ref{l1}, $C_g$ is trivial.  Now, for any $h\in G$,
Ad~$(C_h)\subset C_{{\rm Ad}~(h)}= \{ e \}$.  This implies that
$C_h$ is contained in the kernel of the homomorphism Ad which is the
center of $G$.  Thus, $C_h$ is trivial, hence $G$ is a type $R$ Lie
group. It follows from \cite{Gu} that $G$ has polynomial growth.

We now consider any almost connected group.  Let $G$ be an almost
connected locally compact group.  Then there exists a compact normal
subgroup $K$ such that $G/K$ is a real Lie group.  Let $M$ be the closed
subgroup of $G$ containing $K$ such that $M/K$ is the component of
identity in $G/K$.  Since $G$ is almost connected,
$M$ is a subgroup of finite index.  Now any closed subgroup $M$ is also
unimodular and since $K$ is compact, any closed subgroup of $M/K$ is
unimodular.  Since $M/K$ is a connected Lie group, $M/K$ has polynomial
growth.  Since $K$ is compact, $M$ has polynomial growth.  Since $M$ has finite
index in $G$, $G$ also has polynomial growth.

Converse follows from the facts that closed subgroups of polynomial growth
also have polynomial growth and groups of polynomial growth are unimodular
(cf. \cite{Gu}).
\eo
\end{section}

\begin{tabular}{ll}
Y. Guivarc'h& \hspace*{1cm}  C. R. E. Raja \\
IRMAR& \hspace*{1cm}Stat-Math Unit\\
Campus de Beaulieu & \hspace*{1cm}Indian Statistical Institute\\
Universit\'e de Rennes I & \hspace*{1cm}8th Mile Mysore Road\\
35042 Rennes France & \hspace*{1cm}Bangalore 560059.\\
yves.guivarch@univ-rennes1.fr & \hspace*{1cm}creraja@isibang.ac.in
\end{tabular}


\begin{thebibliography}{Abc}

\footnotesize

\bibitem {Ab} H. Abels, Distal automorphism groups of Lie groups,
J. Reine Angew. Math. 329 (1981), 82--87.

\bibitem {Ba} P. Baldi, Caract\'erisation des groupes de Lie connexes
r\'ecurrents, Ann. Inst. H. Poincar\'e Sect. B (N.S.) 17 (1981), 281--308.

\bibitem {BBE} M. Babillot, P. Bougerol and L. Elie, The random difference
equation $X\sb n=A\sb nX\sb {n-1}+B\sb n$ in the critical case, Ann.
Probab. 25 (1997), 478--493.

\bibitem {BW} U. Baumgartner and Willis, Contraction groups and
scales of automorphisms of totally disconnected locally compact groups,
Israel J. Math. 142 (2004), 221--248.

\bibitem {BQ} Y. Benoist and J.-F. Quint, Mesures stationnaires et ferm\'es 
invariants des espaces homog\`enes, C. R. Acad. Sci. Paris, Ser. I 347 
(2009), 9-13.

\bibitem {BM} S. Blach\`ere and P. Mathieu, Harmonic measures versus
quasiconformal measures for hyperbolic groups, preprint, Arxiv december 2008.

\bibitem {BL} P. Bougerol and J. Lacroix, Products of random matrices with
applications to Schr\"{o}dinger operators. Progress in Probability and
Statistics, 8. Birkh\"{a}user Boston, Inc., Boston, MA, 1985.

\bibitem {Bo} P. Bougerol, Oscillation des produits de matrices al\'eatoires dont
l'exposant de Lyapunov est nul, Lyapunov exponents (Bremen, 1984), 27--36,
Lecture Notes in Math., 1186, Springer, Berlin, 1986.

\bibitem {CKW} D. I. Cartwright, V. A. Kaimanovich and W. Woess, Random walks on
the affine group of local fields and of homogeneous trees. Ann. Inst. Fourier
(Grenoble) 44 (1994), 1243--1288.

\bibitem {CG} J-P. Conze and Y. Guivarc'h, Remarques sur la distalit\'e
dans les espaces vectoriels, C. R. Acad. Sci. Paris S\'er. A  278
(1974), 1083--1086.

\bibitem {DKN} B. Deroin, V. Kleptsyn and A. Navas, On the question of ergodicity
for minimal group action on the circle, Preprint Univ. Rennes 2008.

\bibitem {Du} R. M. Dudley, Random walks on abelian groups, Proc. Amer.
Math. Soc. 13 (1962), 447-450.

\bibitem {EM} A. Eskin and G. Margulis, Recurrence properties of random walks
on finite volume homogeneous manifolds, Random walks and geometry,  431--444,
Walter de Gruyter, Berlin, 2004.

\bibitem {Fo} S. R. Foguel, The ergodic theory of Markov processes. Van Nostrand
Mathematical Studies, No. 21. Van Nostrand Reinhold Co., New York-Toronto,
Ont.-London 1969.

\bibitem{Fur} A. Furman, Random walks on groups and random transformations.
Handbook of dynamical systems, Vol. 1A, 931--1014, North-Holland, Amsterdam, 2002.

\bibitem {Fu} H. Furstenberg, Boundary theory and stochastic processes on
homogeneous spaces. Harmonic analysis on homogeneous spaces (Proc. Sympos. Pure
Math., Vol. XXVI, Williams Coll., Williamstown, Mass., 1972), 193--229. Amer.
Math. Soc., Providence, R.I., 1973.

\bibitem {Fu1} H. Furstenberg, Noncommuting random products. 
Trans. Amer. Math. Soc. 108 (1963), 377--428. 

\bibitem {GS} L. Gallardo and R. Schott, Marches al\'eatoires sur les espaces 
homog\`enes de certains groupes de type rigide.  Conference on Random 
Walks (Kleebach, 1979), pp. 149--170, 4, Ast\'erisque, 74, Soc. Math. France, 
Paris, 1980.

\bibitem {Gr} M. Gromov, Groups of polynomial growth and expanding maps,
Inst. Hautes Études Sci. Publ. Math. No. (1981), 53--73.

\bibitem {Gu} Y. Guivarc'h, Croissance polynomiale et p\'eriodes
des fonctions harmoniques, Bull. Soc. Math. France, 101 (1973), 333--379.

\bibitem {GK} Y. Guivarc'h and M. Keane, Marches al\'eatoires
transitoires et structure des groupes de Lie. Symposia Mathematica, Vol.
XXI (Convegno sulle Misure su Gruppi e su Spazi Vettoriali, Convegno sui
Gruppi e Anelli Ordinati, INDAM, Rome, 1975), 197--217. Academic Press,
London, 1977.

\bibitem {GKR} Y. Guivarc'h, M. Keane and B. Roynette, Marches
al\'eatoires sur les groupes de Lie, Lecture Notes in Mathematics, 624,
Springer-Verlag, Berlin-New York, 1977.

\bibitem {Gu4} Y. Guivarc'h, Quelques propri\'et\'es asymptotiques des produits
de matrices al\'eatoires. Eighth Saint Flour Probability Summer School---1978 (Saint
Flour, 1978), pp. 177--250, Lecture Notes in Math., 774, Springer, Berlin, 1980.

\bibitem {Gu3} Y. Guivarc'h, Propri\'et\'es ergodiques, en mesure infinie,
de certains syst\`emes dynamiques fibr\'es, Ergodic Theory Dynam. Systems
9 (1989), 433--453.

\bibitem {GL} Y. Guivarc'h and Y. LeJan, Sur l'enroulement du flot
g\'eod\'esique, C. R. Acad. Sci. Paris Sér. I Math. 311 (1990), 645--648.

\bibitem {Gu2} Y. Guivarc'h, Marches al\'eatoires sur les groupes,
Development of Mathematics 1950--2000, 577--608, Birkhauser, Basel, 2000.

\bibitem {GuS} Y. Guivarc'h and A. N. Starkov, Orbits of linear group
actions, random walks on homogeneous spaces and toral automorphisms.
Ergodic Theory Dynam. Systems 24 (2004), 767--802.

\bibitem {HR} H. Hennion and B. Roynette, Un th\'eor\`eme de dichotomie pour une
marche al\'eatoire sur un espace homog\`ene. Conference on Random Walks
(Kleebach, 1979) (French), 99--122, 4, Astérisque, 74, Soc. Math. France,
Paris, 1980.

\bibitem{JR}  W. Jaworksi and C. R. E. Raja, The Choquet-Deny
theorem and distal properties of totally disconnected locally compact
groups of polynomial growth, New York J. Math. 13 (2007), 159-174.

\bibitem {Ke} H. Kesten, The Martin boundary of recurrent random walks on
countable groups. 1967 Proc. Fifth Berkeley Sympos. Math. Statist. and
Probability (Berkeley, Calif., 1965/66) Vol. II: Contributions to
Probability Theory, Part 2, pp. 51--74. Univ. California Press, Berkeley,
Calif.

\bibitem {Ke2} H. Kesten, Random difference equations and renewal theory
for products of random matrices, Acta Math. 131 (1973), 207--248.

\bibitem {LeP} E. Le Page, Th\'eor\`emes limites pour les produits de
matrices al\'eatoires, Probability measures on groups (Oberwolfach, 1981),
pp. 258--303, Lecture Notes in Math., 928, Springer, Berlin-New York,
1982.

\bibitem {L} M. Lin, Conservative Markov processes on a topological
space, Israel J. Math. 8 (1970), 165--186.

\bibitem {Lo} V. Losert, On the structure of groups with polynomial
growth. Math. Z. 195 (1987), 109--117.

\bibitem {LM} T. J. Lyons and H. P. McKean, Winding of the plane
Brownian motion, Adv. in Math. 51 (1984), 212--225.

\bibitem {MZ} D. Montgomery and L. Zippin, Topological transformation groups,
Interscience Publishers, New York-London, 1955.

\bibitem {Ma} G. Margulis, Discrete subgroups of semisimple Lie groups,
Ergebnisse der Mathematik und ihrer Grenzgebiete (3) 17. Springer-Verlag, Berlin,
1991.

\bibitem {Mo} C. C. Moore, Ergodicity of flows on homogeneous spaces. Amer. J.
Math. (1966), 154--178.

\bibitem {Pa} A. Parreau, Sous-groupes elliptiques de groupes lin\'eaires sur
un corps valu\'e, J. Lie Theory 13 (2003), 271--278.

\bibitem {Pe} M. Peign\'e, Syst\`emes d'it\'erations de transformations
al\'eatoires: propriet\'es de r\'ecurrence, Preprint, Universit\'e de
Tours, 2008.

\bibitem {Ra2} C. R. E. Raja, On growth, recurrence and the Choquet-Deny
Theorem for $p$-adic Lie groups, Math. Z. 251 (2005), 827--847.

\bibitem {RS} C. R. E. Raja and R. Schott, Recurrent random walks on homogeneous 
spaces of $p$-adic algebraic groups of polynomial growth, Arch. Math. (Basel) 91 
(2008), 379--384. 

\bibitem {Re} D. Revuz, Markov chains. North-Holland Mathematical Library,
Vol. 11. North-Holland Publishing Co., Amsterdam-Oxford; American Elsevier
Publishing Co., Inc., New York, 1975.

\bibitem {Sc} K. Schmidt, Cocycles of ergodic transformation groups.
Macmillan Lectures in Mathematics, Vol. 1. Macmillan Company of India,
Ltd., Delhi, 1977.

\bibitem {Sh} Y. Shalom, Explicit Kazhdan constants for representations of
semisimple and arithmetic groups, Annales de l'institut Fourier, 50
(2000), 833-863

\bibitem {St} G. Stuck, Growth of homogeneous spaces, density of discrete
subgroups and Kazhdan's property (T), Invent. Math. 109 (1992), no. 3, 505--517

\bibitem {Su} D. Sullivan, The density at infinity of a discrete group of
hyperbolic motions, Inst. Hautes Études Sci. Publ. Math. No. 50 (1979), 171--202.

\bibitem {VSC} N. Th. Varopoulos, L. Saloff-Coste and T. Coulhon,
Analysis and geometry on groups, Cambridge Tracts in Mathematics, 100.
Cambridge University Press, Cambridge, 1992.

\bibitem {Wa} S. P. Wang, The Mautner phenomenon for $p$-adic Lie groups,
Math. Z. 185 (1984), 403--412.

\bibitem {Wo} W. Woess,  Random walks on infinite graphs and
groups, Cambridge Tracts in Mathematics, 138. Cambridge University Press,
Cambridge, 2000.



\end{thebibliography}
\end{document}